\documentclass{amsart}

\usepackage{amssymb}
\usepackage{amsmath}

\usepackage{hyperref}
\usepackage{extarrows}

\usepackage{tikz-cd}
\tikzset{
	commutative diagrams/.cd, 
	arrow style=tikz, 
	diagrams={>=stealth}
}

\newtheorem{thm}{Theorem}[section]
\newtheorem{prop}[thm]{Proposition}
\newtheorem{lem}[thm]{Lemma}
\newtheorem{cor}[thm]{Corollary}
\newtheorem{fact}[thm]{Fact}
\newtheorem{thmx}{Theorem}

\theoremstyle{definition}
\newtheorem{definition}[thm]{Definition}
\newtheorem{example}[thm]{Example}

\newtheorem{note}[thm]{Note}

\theoremstyle{remark}

\newtheorem{remark}[thm]{Remark}

\numberwithin{equation}{section}

\newcommand{\R}{\mathbb{R}}  
\newcommand{\C}{\mathbb{C}}  
\newcommand{\A}{\mathbb{A}}  
\newcommand{\Proj}{\mathbb{P}}  
\newcommand{\N}{\mathbb{N}}  
\newcommand{\Z}{\mathbb{Z}}  
\newcommand{\Q}{\mathbb{Q}}  
\newcommand{\F}{\mathbb{F}}  

\newcommand{\squareDiagram}[8]{\begin{equation*}
			\begin{tikzcd}[sep=2.5em, ampersand replacement=\&]
				{#1} \arrow[r, "{#5}"] \arrow[d, swap, "{#6}"] \& {#2} \arrow[d, "{#7}"] \\ 
				{#3} \arrow[r, swap, "{#8}"] \& {#4}
			\end{tikzcd}
		\end{equation*}}

\newcommand{\triangleDiagram}[6]{\begin{equation*}
			\begin{tikzcd}[sep=2.5em, ampersand replacement=\&]
				{#1} \arrow[r, "{#4}"] \arrow[d, swap, "{#5}"] \& {#2}  \\ 
				{#3} \arrow[ur, swap, "{#6}"] \& 
			\end{tikzcd}
		\end{equation*}}

\newcommand{\flippedSquareDiagram}[8]{\begin{equation*}
			\begin{tikzcd}[sep=2.5em, ampersand replacement=\&]
				{#1}  \& {#2} \arrow[l, swap, "{#5}"] \\ 
				{#3} \arrow[u, "{#6}"]  \& {#4} \arrow[l, "{#8}"] \arrow[u, swap, "{#7}"]
			\end{tikzcd}
		\end{equation*}}

\newcommand{\mor}[1]{\ensuremath{\overset{#1}{\longrightarrow}}} 

\DeclareMathOperator{\Ho}{H} 
\DeclareMathOperator{\Poly}{Poly} 
\DeclareMathOperator{\cd}{cd} 
\DeclareMathOperator{\End}{End} 
\DeclareMathOperator{\Hom}{Hom} 
\DeclareMathOperator{\Ind}{Ind} 
\DeclareMathOperator{\Stab}{Stab} 
\DeclareMathOperator{\Id}{Id} 
\DeclareMathOperator{\Aut}{Aut} 

\DeclareMathOperator{\Conf}{Conf} 
\DeclareMathOperator{\Rat}{Rat} 
\DeclareMathOperator{\PConf}{PConf} 
\DeclareMathOperator{\ev}{ev} 

\newcommand{\catname}[1]{{\normalfont\mathbf{#1}}}
\newcommand{\FI}{\catname{FI}}
\newcommand{\VI}{\catname{VI}}
\newcommand{\Arr}{\catname{Arr}}
\newcommand{\Set}{\catname{Set}}
\newcommand{\CCat}{\catname{C}}
\newcommand{\Pos}{\catname{Pos}}
\newcommand{\Vect}{\catname{Vect}}
\newcommand{\Mod}[1]{#1-\catname{Mod}}

\newcommand{\Var}{\catname{Var}}

\renewcommand{\emptyset}{\varnothing}

\begin{document}


\title{Representation stability for families of linear subspace arrangements}


\author{Nir Gadish}
\address{Department of Mathematics, University of Chicago, 
Chicago, Il 29208}
\email{nirg@math.uchicago.edu}





\begin{abstract}
Church-Ellenberg-Farb used the language of $\FI$-modules to prove that the cohomology of certain sequences of hyperplane arrangements with $S_n$-actions satisfies representation stability. Here we lift their results to the level of the arrangements themselves, and define when a collection of \mbox{arrangements is} ``finitely generated". Using this notion we greatly widen the stability \mbox{results to}:
\begin{itemize}
\item General linear subspace arrangements, not necessarily of hyperplanes.
\item A wide class of group actions, replacing $\FI$ by a general category $\CCat$.
\end{itemize}
We show that the cohomology of such collections of arrangements satisfies a strong form of representation stability, with many concrete applications. For example, this implies that their Betti numbers are always given by certain polynomials.

For this purpose we use the theory of representation stability for quite general classes of groups, developed in a companion paper. We apply this theory to get classical cohomological stability of quotients of linear subspace arrangements with coefficients in certain constructible sheaves.
\end{abstract}


\maketitle




\section{Introduction}
A \emph{linear subspace arrangement} is a finite collection $\mathcal{A}$ of linear subspaces in $\C^n$, all containing the origin and possibly of different dimensions. $\mathcal{A}$ determines:

\vspace{1ex} \noindent\textbf{Algebro-geometric data:} An algebraic variety, the complement $M_{\mathcal{A}}=\C^n - \cup \mathcal{A}$.

\vspace{1ex} \noindent\textbf{Combinatorial data:} A partially ordered set $P$ of the intersections of subspaces in $\mathcal{A}$, ordered by reverse inclusion.

\vspace{1ex} \noindent\textbf{Representations:} Cohomology groups $\Ho^i(M_{\mathcal{A}})$ equipped with an action of the group $\Aut(\mathcal{A})$ of $\mathcal{A}$-preserving elements of $\operatorname{GL}_n(\C)$.

\vspace{1ex} \noindent The interaction between these three viewpoints was studied by Arnol'd \cite{Arnol'd}, and later by Goresky-MacPherson \cite{Goresky-MacPherson}, Lehrer-Solomon \cite{Lehrer-Solomon}, and many others.

Many natural arrangements appear in families.
\begin{example} \label{intro:ex1} The \emph{braid arrangement} $B_n = \{ z_i=z_j \}_{1\leq i < j \leq n}$ in $\C^n$, with $\Aut(B_n)=S_n$, the symmetric group on $n$ letters.
\end{example}
\begin{example} \label{intro:ex2}
The arrangement $\displaystyle{C_{(n_1,\ldots, n_m)} = \{ z^{(1)}_{i_1} = \ldots = z^{(m)}_{i_m} \}_{\forall(1\leq j \leq m),\, 1\leq i_j \leq n_j}}$ inside $\displaystyle{\C^{n_1} \times \ldots \times \C^{n_m}}$, with $\displaystyle{\Aut(C_{(n_1,\ldots, n_m)})= S_{n_1}\times \ldots\times S_{n_m}}$.
\end{example}
\begin{example} \label{intro:ex3}
The arrangement $D_n = \{ v_i \neq g(v_j) \}_{1\leq i<j \leq n,\, g\in G}$ in $V^n$, where $G$ is some finite group acting on a complex vector space $V$. Here $\displaystyle{\Aut(D_n)=G^n \rtimes S_n}$.
\end{example}
In the special case of the first example $B_n$, Church, Ellenberg and Farb (\cite{CF} and \cite{CEF}) discovered patterns in the $S_n$-representations $\Ho^i(M_{B_n};\Q)$: their characters can be expressed as a single ``character polynomial" independent of $n$; and their irreducible decompositions stabilize in a precise sense. They named this \emph{representation stability}. The theory of $\FI$-modules, developed in \cite{CEF}, gives a powerful viewpoint that explains this phenomenon as the \emph{finite-generation} of a single object.

\cite{CEF-pointcounts} developed the framework of $\FI$-CHA for discussing families of hyperplane arrangements similar to $B_n$, which captures the sense in which the arrangements themselves are already ``finitely-generated". They then show that the cohomology of such families always forms a finitely-generated $\FI$-module, thus lifting finite-generation to the level of spaces in that case. The approach in \cite{CEF-pointcounts} does not naturally generalize to variations, e.g. to the arrangements in Examples \ref{intro:ex2} and \ref{intro:ex3} above. The purpose of the present paper is to extend these results to include many families of linear subspace arrangements $(A_n)$ (including Examples \ref{intro:ex2} and \ref{intro:ex3}). A major obstacle is that $\Aut(A_n)$ can be quite general. The theory of $\FI$-modules applies only for $\Aut(A_n)=S_n$; it also depends heavily on the specific naming of irreducible representations of $S_n$, which is not available in a more general context. Our goals are three-fold:
\begin{enumerate}
\smallskip\item[\textbf{(I)}] We overcome this obstacle by generalizing the theory of $\FI$-modules to a theory of $\CCat$-modules, where $\CCat$ is any category satisfying certain axioms. When $(A_c)$ is a family of arrangements indexed by $\CCat$, this approach packages $\Ho^i(M_{A_c};\Q)$ into a single object: a \emph{$\CCat$-module}, defined in \S\ref{sec:prelim}. Representation stability then reduces to properties of this one object, namely \emph{finite-generation} and \emph{freeness}.

\smallskip \item[\textbf{(II)}] We extend the stability results of \cite{CEF-pointcounts} to linear subspace arrangements generated by arbitrary linear subspaces, not necessarily hyperplanes.

\smallskip \item[\textbf{(III)}] We \emph{lift} the ``finite-generation" property, which characterizes representation stability of $\FI$-modules, to the level of the spaces themselves.

This project began with a question from Benson Farb who, upon learning that the sequence of complements $M^n - \cup_{i\neq j}\{m_i = m_j\}$ for a manifold $M$ exhibits representation stability when hit with both the cohomology functors $\Ho^i$ and with the homotopy group functors $\pi_i$\footnote{See e.g. Church \cite{Church-manifolds} for $\Ho^i$ and Kupers-Miller \cite{Miller-Kupers} for $\pi_i$.}, asked whether the spaces themselves are finitely-generated in some sense, and the observed finiteness results are a mere shadow of this fact. Here we offer an answer to this question in the case of linear subspace arrangements.
\end{enumerate}

\smallskip The general philosophy behind this work, and it is one that appears repeatedly throughout, is that all representation stability phenomena are the manifestation of ``combinatorial stability" that occurs at the level of certain sets (one might even say that representation stability already exists over $\F_1$). Here this is demonstrated by lifting to the level of various $\CCat$-diagrams of sets and showing that the familiar representation stability follows from this combinatorial stability.

\subsection{Statement of the results}
Let $\CCat$ be a small category. We say that $\CCat$ is \emph{of $\FI$ type}, roughly, if every morphism is a monomorphism and $\CCat$ has pullbacks and push-outs (see \S\ref{step:criterion} below for the precise definition). Categories of $\FI$ type include many natural categories that have recently been studied in the context of representation stability. Among these:
\begin{enumerate}
\item The category $\FI$ itself, of finite sets and injective functions, with automorphism groups $S_n$ for $n\in \N$ (see \cite{CEF}).
\item Finite powers $\FI^m$, with automorphism groups $S_{n_1}\times\ldots \times S_{n_m}$ for  $n_i\in \N$.
\item The category $\VI_k$ of finite dimensional $k$-vector spaces and injective linear maps, with automorphism groups $\operatorname{Gl}_n(k)$ for $n\in\N$ (see Putman-Sam \cite{PS}).
\item The class of categories $\FI_G$ defined by Sam-Snowden \cite{SS-FIG} where $G$ is some group, with automorphism groups $G^n\rtimes S_n$ (see Casto \cite{Ca} for a naturally occurring family of arrangements with these symmetries).
\end{enumerate}

A \emph{$\CCat$-arrangement} $\mathcal{A}_\bullet$ is a functor from $\CCat$ into the category of linear subspace arrangements (see \S\ref{sec:prelim} for the definition of this category). Many infinite families of arrangements can be defined very succinctly and contain a finite amount of information, as illustrated using the following notion. We say that a $\CCat$-arrangement $\mathcal{A}_\bullet$ is \emph{finitely-generated} if there exist finitely many linear subspaces $\{L_i \in \mathcal{A}_{c_i}\}_{i=1}^m$, for $\{c_i\}_{i=1}^m$ objects of $\CCat$, such that for every object $d$ of $\CCat$ the arrangement $\mathcal{A}_d$ is generated by intersecting the images of $\{L_i\}_{i=1}^m$ under all morphisms $c_i\rightarrow d$. Our three examples above each fit into a \mbox{$\CCat$-arrangement} generated by a single subspace: 
\begin{center}
\begin{tabular}{c|c|c|c|c}
Example                                 & Category & Arrangement                 & Group                              & Generating Subspace                               \\ \hline
\ref{intro:ex1} - $B_\bullet$ & $\FI$    & $B_n$                       & $S_n$                              &  				$(z_1=z_2) \subset \C^2$               \\ 
\ref{intro:ex2} - $C_\bullet$ & $\FI^m$  & $C_{(n_1,\ldots,n_m)}$ & $S_{n_1}\times\ldots\times S_{n_m}$ & $(z^{(1)}=\ldots=z^{(m)})\subset \C^m$ \\ 
\ref{intro:ex3} - $D_\bullet$ & $\FI_G$  & $D_n$                       & $G^n\rtimes S_n$                   & $(v_1=v_2) \subset V^2$             \\ 
\end{tabular}
\end{center}

\bigskip Consider a cohomology functor $\Ho^i$. By applying $\Ho^i$ to the complement varieties $M_{\mathcal{A}_d}$ we get a family of abelian groups parameterized by $\CCat$, i.e. a functor from $\CCat$ into abelian groups (or more generally into $\Mod{R}$ for some ring $R$). We call such a functor a \emph{$\CCat$-module}. If $V_\bullet$ is a $\CCat$-module then at every object $d$ of $\CCat$ the module $V_d$ has an action of the group $\Aut_\CCat(d)$, and these representations are related by the morphisms of $\CCat$.

We say that a $\CCat$-module $V_\bullet$ is \emph{finitely-generated} if there exists a finite collection of elements $\{v_i\in V_{c_i}\}_{i=1}^m$ not contained in any proper sub-$\CCat$-module. $V_\bullet$ is \emph{free}\footnote{Free $\FI$-modules are precisely the $\FI\#$-modules presented in \cite{CEF}.} if it is the sum of $\CCat$-modules, each induced from some fixed $\Aut_{\CCat}(d)$-representation for some object $d$ (see \S\ref{step:freeness}). It is the finitely-generated and free $\CCat$-modules that exhibit what we call \emph{representation stability}, i.e. their constituent representations stabilize in a precise sense. The theory of such objects is developed in \cite{Ga} and applied here.

\smallskip Our main result states that the cohomology of a finitely-generated $\CCat$-arrangement is a finitely-generated, free $\CCat$-module. The freeness assertion is the most surprising and consequential part of this statement. We will give concrete applications of this result below. Take $\Ho^*(M)$ to mean singular cohomology with $\Q$-coefficients. For varieties over a field $k$ we also take $\Ho^*(M)$ to be the $\ell$-adic cohomology $\Ho^*_{\acute{e}t}(M_{\bar{k}};\Q_\ell)$ with its \mbox{$\operatorname{Gal}(\bar{k}/k)$-action}.

To avoid pathologies we assume that $\mathcal{A}_\bullet$ respects the structure of $\CCat$, as follows. We say that $\mathcal{A}_\bullet$ is \emph{continuous} if it respects pullbacks in $\CCat$ (see Definition \ref{def:C-arrangements}). $\mathcal{A}_\bullet$ is \emph{normal}, roughly, if no subspace appears later than it could (see Definition \ref{def:normality}). Many natural examples of arrangements satisfy these hypotheses including the ones in Examples \ref{intro:ex1}, \ref{intro:ex2} and \ref{intro:ex3}.

\begin{thmx}
[\textbf{Representation stability of $\CCat$-arrangements}] \label{thm:main_theorem}
Let $\CCat$ be a category of $\FI$ type, and let $\mathcal{A}_\bullet$ be a continuous, normal, and finitely-generated \mbox{$\CCat$-arrangement}. For all $i\geq 0$ the $\CCat$-module of cohomology groups $\Ho^i(\mathcal{M}_{\mathcal{A}})_\bullet$ 
is finitely-generated and free.

When considering $\ell$-adic cohomology, the $\CCat$-module $\Ho^i(\mathcal{M}_{\mathcal{A}})_\bullet$ takes values in the category of $\operatorname{Gal}(\bar{k}/k)$-representations.
\end{thmx}
Theorem \ref{thm:main_theorem} will be proved in \S\ref{sec:main_proof}. We consider two sets of applications:

\smallskip\noindent\textbf{Application I (Representation stability).} For simplicity of exposition we specialize the theory in this introduction to classes of arrangements indexed by \mbox{$m$-tuples} of natural numbers, i.e. $\FI^m$-arrangements and their cohomology. In Example \ref{example:PConf_Conf} we consider an $\FI^m$-family of varieties as follows. Fix $r,k\in \N$, then for every $m$-tuple $(n_1,\ldots, n_m)$ we have the variety $\mathcal{M}^{(n_1,\ldots,n_m)}_{m,k}(\C^r)$ whose points are ordered tuples
$$
[(v^{(1)}_1,\ldots,v^{(1)}_{n_1}),\ldots, (v_1^{(m)},\ldots,v^{(m)}_{n_m})] \in (\C^r)^{n_1}\times\ldots\times (\C^r)^{n_m}
$$ 
where there does not exist any vector $v\in \C^r$ that appears $k$ times in the list ${(v_1^{(j)},\ldots,v^{(j)}_{n_j})}$ for all ${1\leq j\leq m}$. These varieties are the complements of an $\FI^m$-arrangement $\mathcal{A}_{(n_1,\ldots,n_m)}^{m,k}(\C^r)$, generalizing $C_{(n_1,\ldots,n_m)}$ from Example \ref{intro:ex2} above.
\begin{remark}[\textbf{Connection with configuration spaces}]
More geometrically, the varieties $\mathcal{M}^{(n_1,\ldots,n_m)}_{m,k}$ parameterize $m$-tuples of ordered configurations in $\C^r$, where we allow points to collide, subject to the restriction that the configurations cannot all have $k$ points in common including multiplicity. We will discuss special cases of these below.
\end{remark}

In \cite{Ga} we show that for each category $\CCat$ of $\FI$ type we get an algebra of \emph{character polynomials}. These are class functions, simultaneously defined on all groups $\Aut_\CCat(c)$, that uniformly describe the characters of finitely-generated free $\CCat$-modules. In the case $\CCat=\FI^m$ a \emph{character polynomial} is any polynomial in the class functions $X^{(j)}_k$ for $1\leq j \leq m$ and $k\geq 0$, defined simultaneously on all $m$-fold products $\displaystyle{S_{n_1}\times\ldots\times S_{n_m}}$ by
\begin{equation}
X^{(j)}_k(\sigma_1,\ldots,\sigma_m) = \#\text{ of $k$-cycles in } \sigma_j.
\end{equation}

\begin{thmx}[\textbf{Representation stability of $\mathcal{M}^{\bullet}_{m,k}$}] \label{thm:intro_PCOnf}
For every triple of natural numbers $(m,k,r)$, and for each $i\geq 0$, the cohomology
$$
\Ho^i(\mathcal{M}^\bullet_{m,k}(\C^r))
$$
forms a finitely-generated, free $\FI^m$-module of multi-degree $\lfloor \frac{i}{r}\rfloor k(1,\ldots,1)$. 

In particular, there exists a single $\FI^m$-character polynomial 
$$
P_i\in \Q[X^{(j)}_d \mid 1\leq j \leq m, d\geq 0]
$$
of multi-degree $\lfloor \frac{i}{r}\rfloor k(1,\ldots,1)$ such that for every $(n_1,\ldots,n_m)\in \N^m$ there is an equality of class functions
\begin{equation}
\chi_{\Ho^i(\mathcal{M}^{(n_1,\ldots,n_m)}_{m,k}(\C^r))} = P_i
\end{equation}
The multi-degree of $X^{(j)}_d$ is defined to be $d \cdot \bar{e}^{(j)}$.
\end{thmx}
\begin{example} When $m=2$ and $k=r=1$ the varieties $\mathcal{M}^{(n_1,n_2)}_{2,1}(\C)$ are covers of the space of rational maps $\Rat^n_*(\C)$ studied by Segal (see Example \ref{ex:rational_m} below). In this case, the character polynomials $P_1$ and $P_2$ are $\displaystyle{\chi_{\Ho^1(\mathcal{M}^{(n_1,n_2)}_{2,1}(\C))}=X^{(1)}_1\cdot X^{(2)}_1}$ of multi-degree $(1,1)$ and
\begin{eqnarray*}
\chi_{\Ho^1(\mathcal{M}^{(n_1,n_2)}_{2,1}(\C))} = X^{(1)}_1\left( \binom{X^{(2)}_1}{2} -X_2^{(2)}  \right) + X^{(2)}_1\left( \binom{X^{(1)}_1}{2} -X_2^{(1)}  \right)\\ + 2\binom{X^{(1)}_1}{2}\binom{X^{(2)}_1}{2}-2X^{(1)}_2X^{(2)}_2
\end{eqnarray*}
of multi-degree $(2,2)$, both independent from $(n_1,n_2)$.
\end{example}

Other special cases to which Theorem \ref{thm:intro_PCOnf} applies include:
\begin{enumerate}
\item The braid arrangements, i.e. the classifying spaces of Artin's braid groups, discussed in Example \ref{ex:braid} below.
\item Spaces of configurations in $\C^r$, discussed in Example \ref{ex:Conf_r} below.
\item The $k$-equals arrangements, related to incomputability problems and classifying homotopy links, discussed in Example \ref{ex:k_equals} below.
\item Covers of the spaces of based holomorphic maps $\Proj^1\mor{}\Proj^m$, discussed in Example \ref{ex:rational_m} below.
\end{enumerate}
Theorem \ref{thm:intro_PCOnf} shows that all of these examples exhibit representation stability with explicit stable ranges. Moreover, we get information regarding the Betti numbers of all these varieties. Applying Theorem \ref{thm:intro_PCOnf} to $\sigma = id$ gives the following.
\begin{cor}[\textbf{Polynomial Betti numbers for $\mathcal{M}^\bullet_{m,k}$}]
For every $i\geq 0$ there exists a polynomial $p_i\in\Q[t_1,\ldots,t_m]$ of multi-degree $\lfloor\frac{i}{r}\rfloor k(1,\ldots,1)$ such that
$$
\dim_{\Q} \Ho^i(\mathcal{M}^{(n_1,\ldots,n_m)}_{m,k}(\C^r)) = p_i(n_1,\ldots,n_m)
$$
for all $(n_1,\ldots,n_m)\in \N^m$.
\end{cor}

\bigskip\noindent\textbf{Application II (Classical cohomological stability).} The quotient spaces  $\displaystyle{X_{A_d}:= M_{A_d}/\Aut_\CCat(d)}$ come up in multiple contexts:
\begin{enumerate}
\item For $B_n$ from Example \ref{intro:ex1}, $X_{B_n}$ is the space $\Poly_n(\C)$ of degree-$n$ square-free polynomials, studied by Arnol'd \cite{Arnol'd}.
\item For $C_{(n,\ldots,n)}$ from Example \ref{intro:ex2}, $X_{C_{(n,\ldots,n)}}$ is the space $\operatorname{Hol}^n_*(\Proj^1,\Proj^m)$ of based degree-$n$ holomorphic maps, studied by Segal \cite{Se}.
\item For $\mathcal{M}^{(n,\ldots,n)}_{m,k}(\C)$ discussed above, the quotient $X^{n}_{m,k}$ is the space $\Poly^n_{m,k}(\C)$ of $m$-tuples of polynomials with restrictions on root coincidences, introduced by Farb-Wolfson \cite{FW}. These generalize the two previous examples.
\end{enumerate}
The cohomology of a quotient of some variety $M$ by a finite group $G$ is given by transfer
\begin{equation}
\Ho^i(M/G;\Q) = \Ho^i(M;\Q)^G.
\end{equation}
Thus Theorem \ref{thm:main_theorem} applied to the trivial subrepresentation gives a classical cohomological stability statement.
\begin{thmx}[\textbf{Cohomological stability for arrangement quotients}]\label{thm:intro-cohomological}
Suppose that $\CCat$ and $\mathcal{A}_\bullet$ satisfy the hypotheses of Theorem \ref{thm:main_theorem} and that $|\Aut_\CCat(d)|<\infty$ for every object $d$. Then for both singular and \'{e}tale cohomology, the groups $\Ho^i(X_{\mathcal{A}_d})$ stabilize in the following sense: if any morphisms $c\rightarrow d$ exist in $\CCat$ then there is a canonical injection
$$
\Ho^i(X_{\mathcal{A}_c}) \hookrightarrow \Ho^i(X_{\mathcal{A}_d})
$$ 
and these maps become isomorphisms when $c$ is sufficiently large relative to $i$ (an explicit stable range is given in \S\ref{sec:cohomological_stability}).
\end{thmx}
In some special cases, Theorem \ref{thm:intro-cohomological} was previously proved for integral cohomology using clever but ad-hoc techniques, see e.g. \cite{Se}. Applying the theorem to $\mathcal{M}^{(n,\ldots,n)}_{m,k}$ we get a new proof of the cohomological stability proved in \cite{FW} for the rational cohomology of $\Poly^n_{m,k}(\C)$.

However, Theorem \ref{thm:intro-cohomological} considers only the trivial subrepresentation of $\Ho^i(\mathcal{M}_{\mathcal{A}_d})$. This is a very special case of the following more general version that considers the entire representation. Every $\CCat$-module $N_\bullet$ induces a natural constructible sheaf $\widetilde{N}_c$ on the quotient  $X_{\mathcal{A}_d}$ whose stalk above the orbit $[x]$ satisfies
$$
(\widetilde{N}_d)_{[x]} \cong N_d^{\Stab(x)}
$$
where $\Stab(x)$ is the $\Aut_\CCat(d)$-stabilizer of $x\in \mathcal{M}_{\mathcal{A}_d}$ (see \S\ref{sec:cohomological_stability} for the construction).
\begin{thmx}[\textbf{Twisted stability for arrangement quotients}]\label{thm:intro-twisted-cohomological}
Suppose that $\CCat$ and $\mathcal{A}_\bullet$ satisfy the hypotheses of Theorem \ref{thm:intro-cohomological}. Let $N_\bullet$ be a finitely-generated, free $\CCat$-module. Then the sheaf cohomology groups $\Ho^i(X_{\mathcal{A}_d};\widetilde{N}_d)$ stabilize in the sense of Theorem \ref{thm:intro-cohomological}.

When working with $\ell$-adic cohomology one needs $N_\bullet$ to take values in continuous $\Q_\ell$-modules.
\end{thmx}
By considering the trivial $\CCat$-module $N_\bullet\equiv\Q$ we recover Theorem \ref{thm:intro-cohomological}. The proof of Theorem \ref{thm:intro-twisted-cohomological} will be presented in \S\ref{sec:cohomological_stability}.

Lastly, through the Grothendieck-Lefschetz fixed-point theorem, the twisted stability result in Theorem \ref{thm:intro-twisted-cohomological} has implications for arithmetic statistics of varieties over finite fields. This direction will be developed further in a sequel \cite{CG}.

\subsection{Acknowledgments}
I would like to thank Kevin Casto, Weiyan Chen and Jonathan Rubin for many hours of helpful conversation. Also, I am deeply grateful to my advisor Benson Farb for presenting me with the problem, and to Farb and my second advisor Jesse Wolfson who, through their many helpful comments and suggestions, made this paper take the shape it now has.

\section{Preliminaries} \label{sec:prelim}
This section will introduce the necessary terminology and categories in which we will be working. We start with the fundamental object of interest, namely linear subspace arrangements.
\begin{definition}[\textbf{Linear arrangements}]
The category of \emph{linear subspace arrangements} over a field $k$, denoted by $\Arr_k$, consists of pairs $A=(V,L)$ where $V$ is a finite dimensional vector space over $k$ and $L$ is a finite set of linear subspaces of $V$, all containing the origin\footnote{This is often called \emph{a central arrangement}. All arrangements here are assumed to be central.}, such that $L$ is closed under intersections. A \emph{morphism} $(V_1, L_1) \mor{f} (V_2,L_2)$ is a surjective linear map $V(f):V_2\mor{}V_1$ such that for every subspace $W\in L_1$ the preimage $V(f)^{-1}(W)$ belongs to $L_2$.
\end{definition}
As the definition suggests, there are two natural functors from the category of arrangements:
\begin{definition}[\textbf{Underlying vector space}]
 The \emph{underlying vector spaces} functor $V: \Arr_k^{op} \mor{} \Vect_k$ is defined by sending an arrangement $(V,L)$ to the vector space $V$. Morphisms of arrangements are defined as being (contravariant) linear maps between the vector spaces, and these define the action of $V$ on morphisms: $f\mapsto V(f)$.
\end{definition}
and secondly,
\begin{definition}[\textbf{Intersection semilattice}] The \emph{intersection semilattice} functor $L: \Arr_k \mor{}\Pos$ is defined by sending an arrangement $(V,L)$ to the ranked poset $(L,\cd)$ ordered by reverse inclusion of subspaces, where $\cd$ is the codimension function: $\cd(W) = \dim_k(V)-\dim_k(W)$ (see \ref{def:category_pos} later for the definition of the category $\Pos$ of ranked posets). A morphism $(V_1,L_1)\mor{f}(V_2,L_2)$ defines a set function $L_1\mor{L(f)}L_2$ by $W\mapsto V(f)^{-1}(W)$, which preserves inclusions and respects intersection.
\end{definition}

An arrangement naturally gives rise to an algebraic variety:
\begin{definition}[\textbf{The complement of an arrangement}]
The \emph{complement functor} $\displaystyle{\mathcal{M}:\Arr_k^{op} \mor{} \Var_k}$ from arrangements to the category of algebraic varieties over $k$ (or when $k=\C$, to complex manifolds) is the contravariant functor that sends an arrangement $(V, L)$ to the $k$-variety $V - \cup L$. The morphism of arrangements $(V_1,L_1)\mor{f}(V_2,L_2)$ induces a map $V_2-\cup L_2 \mor{V(f)} V_1 - \cup L_1$ by restriction.
\end{definition}
Note that the preimage of $\cup L_1$ is contained inside $\cup L_2$, and thus the restriction is well-defined.

We are concerned with families of arrangements, and their complements, indexed by some category $\CCat$. Formally this is given by a diagram of arrangements, i.e. a functor. Throughout this paper we denote a covariant (resp. contravariant) functor $\mathcal{F}: X \mor{} Y$ by $\mathcal{F}_\bullet$ (resp. $\mathcal{F}^\bullet$).
\begin{definition} [\textbf{$\CCat$-arrangements}] \label{def:C-arrangements}
Let $\CCat$ be a category. A $\CCat$-arrangement (over $k$) is a functor $\mathcal{A}: \CCat \mor{} \Arr_k$. We denote the compositions $V\circ \mathcal{A}$, $L\circ \mathcal{A}$ and $\mathcal{M}\circ \mathcal{A}$ by $V_\mathcal{A}^\bullet$, $L^\mathcal{A}_\bullet$ and $\mathcal{M} _\mathcal{A}^\bullet$ respectively. Note that the associations ${c\mapsto \mathcal{A}_c\mapsto V(\mathcal{A}_c), \mathcal{M}(\mathcal{A}_c)}$ are contravariant and are therefore denoted with an upper index, i.e. $\mathcal{A}_c = (V_\mathcal{A}^c, L^\mathcal{A}_c)$.

If $\mathcal{A}$ is a $\CCat$-arrangement, we say that the underlying diagram of vector spaces $V_{\mathcal{A}}$ is \emph{continuous} if it takes pullback diagrams in $\CCat$ to push-out diagrams in $\Vect_k$. When this is the case, we also say that $\mathcal{A}_\bullet$ itself is \emph{continuous}.
\end{definition}
We will often omit the superscript and subscript of $\mathcal{A}$ from $L^{\mathcal{A}}$ and $V_{\mathcal{A}}$ when there is no ambiguity as to which arrangement is involved.

By applying a cohomology functor $\Ho^i$ to the $\CCat^{op}$-variety $\mathcal{M}_\mathcal{A}$ we get a $\CCat$-module, and these modules are the subject of our Main Theorem \ref{thm:main_theorem}. They are defined as follows. Note that since $\mathcal{M}$ and $\Ho^i$ are both contravariant functors, their composition is covariant.
\begin{definition}[\textbf{$\CCat$-modules and finite-generation}]
Let $R$ be some ring and $\Mod{R}$ its category of modules. A covariant functor $M: \CCat\mor{} \Mod{
R}$ is called \emph{a $\CCat$-module over $R$}, or more implicitly \emph{a $\CCat$-module}. The category of all $\CCat$-modules naturally forms an abelian category with kernels and cokernels computed pointwise. In particular, there are natural notions of $\CCat$-submodules and quotients, defined by pointwise injections and surjections respectively.

We say that a $\CCat$-module $M_\bullet$ is \emph{generated} by a set $X=\{x_\alpha \in M_{c_\alpha} \}_{\alpha\in A}$ if there does not exist any proper $\CCat$-submodule of $M_\bullet$ that contains $X$. Equivalently, $M_\bullet$ is generated by $X$ if for every object $c$ of $\CCat$ the $R$-module $M_c$ is generated by the images of $X$ under maps induced from morphisms of $\CCat$.

If $M_\bullet$ is generated by a finite set, we say that it is \emph{finitely-generated}.
\end{definition}

\section{Steps towards representation stability} \label{sec:main_proof}
We now set out to prove Theorem \ref{thm:main_theorem}. As discovered by Goresky-MacPherson, the cohomology of an arrangement complement is determined by the combinatorial data encoded in its partially ordered set of subspaces. We therefore start the proof by setting up the terminology involving these objects and the notion of their combinatorial stability.

\subsection{Step 1 - Ranked posets and combinatorial stability} \label{step:ranked_posets}
\begin{definition}[\textbf{Ranked posets}] \label{def:category_pos}
The category of \emph{finite ranked posets}, denoted by $\Pos$, is described as follows. The objects are pairs $(P,r)$ where $P$ is a finite partially ordered set and $r$ is a function $r: P\mor{}\Z$, called the \emph{rank function}, that is strictly increasing. A morphism $(P_1, r_1) \mor{f} (P_2,r_2)$ is a set function $f:P_1\mor{}P_2$ that preserves both ordering and rank, i.e. $x\leq_1 y \implies f(x)\leq_2 f(y)$ and $r_1(x) = r_2(f(x))$ for all $x,y\in P_1$.
\end{definition}
In our context the rank will often be determined by codimension.

\bigskip Recall that the homology of a poset $P$ is defined to be the homology of its nerve, i.e. the simplicial set $\Delta(P)$ whose $n$-simplices are order-chains ${x_0\leq x_1\leq \ldots \leq x_n}$ in $P$. Taking the nerve $P\mapsto \Delta(P)$ and the poset homology are functors. The following concept explicitly appears in Goresky-MacPherson's formula for the cohomology of the complement of a linear subspace arrangement. 
\begin{definition}[\textbf{Whitney homology}]
The \emph{Whitney homology} functors of a ranked poset $(P,r)$ are defined by
\begin{equation}
WH_n(P,r) := \bigoplus_{x\in P^n} \widetilde{\Ho}_{n-2}(P^{<x})
\end{equation}
where $\widetilde{\Ho}$ stands for reduced homology, $n$ is an integer, $P^n$ is the subposet of elements with rank $n$, and $P^{<x} = \{y\in P \mid y<x\}$ with its induced ordering.

A morphism $(P_1,r_2)\mor{f}(P_2,r_2)$ sends every subposet $P_1^{<x}$ to $P_2^{<f(x)}$ and thus induces homomorphisms $\widetilde{\Ho}_{n-2}(P_1^{<x}) \mor{f_*} \widetilde{\Ho}_{n-2}(P_2^{<f(x)})$. The direct sum of these homomorphisms is the induced homomorphism $WH_n(P_1,r_1)\mor{f_*}WH_n(P_2,r_2)$.
\end{definition}
\begin{remark}
When $(P,r)$ is a geometric lattice, the Whitney homology groups defined above coincide with the homology of the complex $\CCat_\bullet(\Delta(P))$ with the truncated differential
$$
\partial^W_n = \sum_{i=0}^{n-1}(-1)^i\delta_i
$$
at level $n$. See \cite{Bj} for a discussion on the motivation for this definition and a proof of this claim.
\end{remark}

We now consider families of posets. Let $\CCat$ be some indexing category.
\begin{definition}
A $\CCat$-poset is a functor $P_\bullet: \CCat \mor{} \Pos$.
\end{definition}

Combinatorial stability of such families of posets is defined by the following two properties.
\begin{definition}[\textbf{Finite generation}]
We say that a $\CCat$-poset $P_\bullet$ is \emph{finitely-generated} if for every rank $n\in \Z$ there exist finitely many elements $\{x_i\in P_{c_i}^{n}\}_{i=1}^{k}$ whose orbits under $\CCat$ contain $P^n_d$ for every object $d$.
\end{definition}
\begin{note}
If we change the rank function by composing it with an injective order preserving function $\Z\mor{}\Z$, the notion of being finitely-generated remains unchanged.
\end{note}

\begin{definition}[\textbf{Downward stability}]
We say that a $\CCat$-poset $P_\bullet$ is \emph{downward stable} if for every morphism $c\mor{f}d$ in $\CCat$ and an element $x\in P_c$ the induced poset map $P_c^{<x}\mor{f_*}P_d^{<f(x)}$ is an isomorphism.
\end{definition}
We can now phrase our notion of stability for $\CCat$-posets.
\begin{definition}[\textbf{Combinatorial stability}]
The $\CCat$-poset is said to exhibit \emph{combinatorial stability} if it is both finitely-generated and downward stable.
\end{definition}

As the following theorem shows, combinatorial stability at the level of $\CCat$-posets implies more familiar stability phenomena that occur in the context of representation stability: recall that an $\FI$-module is representation stable in the sense of \cite{CF} if and only if it is finitely-generated (see \cite{CEF}).
\begin{thm}[\textbf{Finitely-generated $\CCat$-poset homology}] \label{thm:combinatorially_stable_implies_fg}
If a $\CCat$-poset $(P,r)_\bullet$ is combinatorially stable, then its Whitney homology $WH_n(P,r)_\bullet$ is a finitely-generated $\CCat$-module for all $n$.
\end{thm}
\begin{proof}
Let $n$ be a natural number. Since $P_\bullet$ is finitely-generated we can find a finite list of elements $x_i\in P_{c_i}^n$ for $i=1,\ldots,l$ whose $\CCat$-orbits contain all rank $n$ elements of $P_\bullet$. Each of the Whitney homology groups $WH_n(P_{c_i},r_{c_i})$ is finitely-generated, so it will suffice to show that their $\CCat$-orbits span all other Whitney homology groups.

Let $d$ be any object of $\CCat$ and $y\in P_d^n$. Suffice it to show that $\widetilde{\Ho}_{n-2}(P_d^{<y})$ is contained in the $\CCat$-orbits of the above groups. By our choice of $x_1,\ldots, x_l$, there exists some $1\leq i \leq l$ and a morphism $c_i \mor{f} d$ such that $f_*(x_i)=y$. By the downward stability assumption, the induced map
$$
P_{c_i}^{<x_i}\mor{f_*} P_d^{<y}
$$
is an isomorphism, and therefore the induced homomorphism on homology is also an isomorphism. In particular it is surjective.
\end{proof}

\subsection{Step 2 - The cohomology of an arrangement complement} \label{step:cohomology_of_arrangements}
Goresky-MacPherson used Stratified Morse Theory to give a formula for the cohomology groups of real and complex linear subspace arrangement complements in terms of the associated intersection semilattice (see \cite{Goresky-MacPherson}). Later, Bj\"{o}rner-Ekedahl compute the $\ell$-adic cohomology of the complement of a linear subspace arrangement defined over some arbitrary field $k$ (see \cite{BE}). Their formula coincides with the Goresky-MacPherson result for the case $k=\C$. The cohomology groups are given as follows.
\begin{thm}[\cite{BE}, Theorem 4.9]
Let $A = (V,L)$ be a subspace arrangement and $\mathcal{M}_A$ its complement. If $\ell\neq char(k)$ is a prime number then the $\ell$-adic cohomology of $\mathcal{M}_A$ is given by
\begin{eqnarray} \label{eq:Bjorner_Ekedahl}
\Ho^i_{\acute{e}t}(\mathcal{M}_A; \Q_\ell) &\cong& \bigoplus_{x\in L^A}\widetilde{\Ho}_{2\cd(x)-i-2}(\Delta(L^{<x}))\otimes \Q_\ell(\cd(x))\\ \label{eq:Bjorner-Whitney_sum}
& = & \bigoplus_{n\geq 0} WH_{2n-i}(L, 2\cd - i)\otimes \Q_\ell(n)
\end{eqnarray}
where $\cd(x)$ is the codimension of the subspace $x$ in $V$, and $L^{<x}$ is the subposet of spaces in $L$ that contain $x$. The term $\Q_\ell(n)$ is the $\ell$-adic cyclotomic character (see \cite{tate-twist}).

\end{thm}
The key feature of Isomorphism \eqref{eq:Bjorner_Ekedahl} is that it is natural, i.e. we can read off pullback maps between cohomology groups from the poset maps and the induced homomorphisms on Whitney homology. This can be seem e.g. by applying Poincar\'{e} duality to the spectral sequence introduced by Petersen in \cite{Petersen}.

Implicit in Equation \eqref{eq:Bjorner-Whitney_sum} is that the direct sum of Whitney homologies is finite. This observation is essential to establishing finite generation of $\CCat$-modules that occur as cohomology of complements of $\CCat$-arrangements. Indeed,
\begin{lem}\label{lem:Deligne_bounds}
The direct sum decomposition for $\Ho^i_{\acute{e}t}(\mathcal{M}_A; \Q_\ell)$ given in \eqref{eq:Bjorner_Ekedahl} and \eqref{eq:Bjorner-Whitney_sum} includes contributions only from subspaces $x\in L$ of codimension ${\frac{i}{2}\leq \cd(x) \leq i}$. Equivalently, the only Whitney homology groups that contribute to the direct sum are the ones whose index $n$ satisfies $\frac{i}{2}\leq n \leq i$.
\end{lem}
\begin{proof}
By Deligne's bounds \cite{De} the weights that occur in $\Ho^i_{\acute{e}t}(X;\Q_\ell)$ are bounded between $i$ and $2i$. The lemma now follows from the fact that the weight of the $n$-th summand of \eqref{eq:Bjorner_Ekedahl} is $2n$.

Alternatively, we can see this directly by elementary means. For the lower bound notice that if $x\in L$ has codimension smaller than $\frac{i}{2}$, then the direct summand corresponding to $x$ in \eqref{eq:Bjorner_Ekedahl} is
$$
\widetilde{\Ho}_{2\cd(x)-i-2}(\Delta(L^{<x}))\otimes \Q_\ell(\cd(x))
$$
where $2\cd(x)-i-2 \leq -2$. Since reduced homology groups are zero below degree $-1$, this summand is zero
\footnote{Note that there could be a contribution in degree $-1$ since, by convention, $\widetilde{\Ho}_{-1}(\emptyset;\Z)=\Z$. This term occurs precisely when $x$ is minimal in $L$, or equivalently not contained in any other subspace in $L$.}.

For the upper bound the claim will follow if we show that $\Delta(L^{<x})$ has non-degenerate $(2\cd(x)-i-2)$-simplices only when $\cd(x)\leq i$. Suppose
$$
x_0<x_1<\ldots< x_n  \, (< x)
$$
is a non-degenerate simplex in $\Delta(L^{<x})$. Then the strict monotonicity of codimension gives
$$
1\leq \cd(x_0) < \cd(x_1) < \ldots < \cd(x_n) < \cd(x)
$$
and since these are all integers, $n+1 < \cd(x)$. Thus the existence of a non-degenerate $(2\cd(x)-i-2)$-simplex implies that $(2\cd(x)-i-2)+1 < \cd(x)$, or equivalently $\cd(x) < i+1$.
\end{proof}

\begin{cor}[\textbf{Stable semilattices imply controlled cohomology}]\label{cor:stable_lattice_implies_fg}
Let $\mathcal{A}: \CCat \mor{} \Arr_k$ be a $\CCat$-arrangement. If the intersection semilattice $L^\mathcal{A}_\bullet$ is combinatorially stable (i.e. finitely-generated and downward stable), then for each $i\geq 0$ the $\CCat$-module $\Ho^i_{\acute{e}t}({\mathcal{M}_\mathcal{A}}; \Q_\ell)_\bullet$ is finitely-generated.
\end{cor}
\begin{proof}
We have seen in Theorem \ref{thm:combinatorially_stable_implies_fg} that a combinatorially stable $\CCat$-poset gives rise to finitely-generated Whitney homology groups (in every degree). Since $\Ho^i(\mathcal{M}_\mathcal{A})$ is naturally isomorphic to a \textbf{finite} direct sum of such homology groups, the resulting cohomology $\CCat$-module is finitely-generated.
\end{proof}

\subsection{Step 3 - Criterion for stability} \label{step:criterion}
This subsection discusses properties of $\CCat$-arrangements which, together with some structural assumptions on the category $\CCat$, will ensure that the associated intersection semilattice will be combinatorially stable. First we define a notion of finite-generation for $\CCat$-arrangements. It is this property of an arrangement that ensures the finite-generation of the associated intersection semilattice. Downward stability poses more of a challenge, and it will lead us to the notion of a normal $\CCat$-arrangement.
\begin{definition}[\textbf{Finitely-Generated $\CCat$-arrangements}]
A $\CCat$-arrangement $\mathcal{A}_\bullet$ is said to be \emph{generated by the set of subspaces} $\left\{x_\alpha \subset V^{c_\alpha}\right\}_{\alpha\in A}$ if for every object $d$ of $\CCat$ and every subspace $y\in L^\mathcal{A}_d$ there exists a finite list of morphisms $c_{\alpha_i}\mor{f_i}d$ where $1\leq i \leq l$ such that
\begin{equation}
y = L(f_1)x_{\alpha_1}\cap \ldots \cap L(f_l)x_{\alpha_l}.
\end{equation}
Equivalently, if $\mathcal{A}_\bullet$ is the least $\CCat$-arrangement that contains all of the subspaces $\{x_\alpha\}_{\alpha\in A}$ among its chosen subspaces.

When this is the case, we say that $\mathcal{A}_\bullet$ is \emph{generated in degrees} $\{c_\alpha\}_{\alpha\in A}$. The \mbox{$\CCat$-arrangement} $\mathcal{A}_\bullet$ is \emph{finitely-generated} if it is generated by some finite set of subspaces.
\end{definition}
\begin{example}
The braid $\FI$-arrangement over $k$ is generated in degree $2$ by a single subspace: $\{z_1=z_2\} \subset k^2$. See example \ref{ex:braid} for an elaboration.
\end{example}

The following notation will prove useful.
\begin{definition}\label{def:ordering_objects}
Let $\CCat$ be a category. If $c$ and $d$ are two objects, we say that $c\leq d$ if $\Hom_\CCat(c,d)\neq \emptyset$. Moreover, we say that $c < d$ if $c\leq d$ and $d \not\leq c$.
\end{definition}

Downward stability turns out to be related to a notion of saturation of a $\CCat$-arrangement. We make this connection precise using the following definition.
\begin{definition}[\textbf{Normality and primitive subspaces}] \label{def:normality}
Let $\mathcal{A}_\bullet = (V^\bullet,L_\bullet)$ be a $\CCat$-arrangement.
\begin{itemize}
\item $\mathcal{A}_\bullet$ is \emph{normal} if for every morphism $c\mor{f}d$, every subspace $x\in L^\mathcal{A}_d$ that contains $\ker V(f) $ is in the image of $L(f): L^\mathcal{A}_c\mor{}L^\mathcal{A}_d$. Equivalently, the direct image $V(f)x \subseteq V^c$ is already in $L^\mathcal{A}_c$.

\item A subspace $x\in L^\mathcal{A}_d$ is \emph{primitive} if it does not contain the kernel of any linear map induced by a morphism $c\mor{}d$ where $c<d$. We define the \emph{degree} of $x$ to be the object $d$ and denote $\deg(x)=d$.
\end{itemize}
Normality and primitivity are well-behaved since we are assuming that all the subspaces in $L^\mathcal{A}$ contain the origin.
\end{definition}
\begin{example}
To illustrate the meaning of (non)normality, consider the following example. Let $\CCat = \{ 0\mor{f}1 \}$ be a category of two objects with a single morphism between them. Define a diagram of vector spaces by $V^0=V^1 = \C$ with ${V(f)=\Id: V^1 \mor{} V^0}$ and construct a $\CCat$-arrangement $\mathcal{A}_\bullet = (L_\bullet, V^\bullet)$ by choosing
$$
L_0 = \emptyset, \; L_1 = \{0\}.
$$
The arrangement thus constructed is not normal, as $0\in L_1$ is not in the image of $L(f)$ even though it could be there (and perhaps ``morally should" be there). Had we chosen $L_0=L_1$ we would have defined a normal arrangement, since then every element of $L_1$ that \emph{could} lie in the image of $L(f)$ indeed appears there.
\end{example}

As previously declared, normality guarantees downward stability. 
\begin{lem} \label{lem:normal_is_stable}
The intersection semilattice $L^\mathcal{A}_\bullet$ of a normal $\CCat$-arrangement $\mathcal{A}$ is downward stable.
\end{lem}
\begin{proof}
Throughout this proof we will denote the poset $L^\mathcal{A}_c$ by $L_c$. Suppose ${c\mor{f}d}$ and $x\in L_c$. We will construct an inverse to the induced poset morphism $${L_c^{<x}\mor{L(f)} L_d^{<L(f)x}}.$$ By definition of the order on $L_d$, every $y\in L_d^{<L(f)x}$ satisfies $L(f)(x)\subseteq y$ and thus
$$
\ker V(f) \subseteq V(f)^{-1}(x) \subseteq y,
$$
so by normality $V(f)y$ is a subspace in $L_c$. Since $V(f)$ is surjective we have ${x\subseteq V(f)y}$, so $V(f)y\in L_c^{<x}$. The inverse map to $L(f) = V(f)^{-1}$ is therefore the direct image $y \mapsto V(f)y$.
\end{proof}

Next we show that under structural assumptions on $\CCat$ the properties defined above indeed ensure the combinatorial stability of the associated intersection semilattice.
\begin{definition}[\textbf{Weakly filtering categories}]
A category $\CCat$ is \emph{weakly filtering} if for every pair of objects $c_1$ and $c_2$ there exists a finite collection of objects $d_1,\ldots, d_k$ and morphisms $c_i\mor{f_{ij}}d_j$ where $i=1,2$ and $1\leq j \leq k$, such that every pair of morphisms $c_i\mor{g_i}e$ for $i=1,2$ factors through one of the $d_j$'s. Explicitly, for every pair of morphisms $c_i\mor{g_i}e$ there exists some $1\leq j\leq k$ and a morphism $d_j\mor{g}e$ such that $g\circ f_i = g_i$ for $i=1,2$.

This property is called \emph{property (F)} by Sam-Snowden in \cite{SS}. 
\end{definition}
\begin{lem} \label{lem:filtering_is_fg}
If the category $\CCat$ is weakly filtering, then the intersection semilattice of every finitely-generated $\CCat$-arrangement is a finitely-generated $\CCat$-poset.
\end{lem}
\begin{proof}
Fix a finite set of generators $X = \left\{ x_\alpha \subset V^{c_\alpha} \right\}_{\alpha\in A}$ and a codimension $n$. For every object $d$ and a subspace $y\in L_e$ of codimension $n$, there exists a finite list of morphisms $c_{\alpha_i}\mor{g_i}d$ such that
$$
y = L(g_1)x_{\alpha_1} \cap \ldots \cap L(g_l)x_{\alpha_l}.
$$
But since $y$ is of codimension $n$, it can be written as the intersection of no more than $n$ subspaces. Thus without loss of generality we can assume that $l\leq n$. We will show that all such intersections are in the $\CCat$-orbits of $\bigcup_{d\in D_n} L_d$ for a fixed finite collection $D_n$ of objects in $\CCat$. Since every poset $L_d$ is finite, the union $\bigcup_{d\in D_n} L_d$ has finitely many elements, and this will complete the proof.

We prove this by induction on $n$. For $n=1$ it suffices to take the finite set $D_1=\{ c_\alpha \}_{\alpha\in A}$. For the induction step, suppose $D_{n-1}$ is already defined. For every $c_1\in D_{n-1}$ and $c_2\in D_1$ find a finite collection of objects $d_1,\ldots, d_k$ through which every pair of morphisms factors. The set $D_n$ will be defined to be the union of these finite lists as $c_1$ and $c_2$ range over the finite sets $D_{n-1}$ and $D_1$ respectively.

We need to show that $D_n$ satisfies the desired property. Suppose $y\in L_e$ is given by
$$
y = L(g_1)x_{\alpha_1} \cap \ldots \cap L(g_l)x_{\alpha_l}
$$
as above with $l\leq n$. By repeating the last term if necessary, we can assume that $l=n$. By the choice of $D_{n-1}$ there exists some $c_1\in D_{n-1}$ and a morphism $c_1\mor{h_1}e$ such that the $(n-1)$-fold intersection $L(g_1)x_{\alpha_1} \cap \ldots \cap L(g_{n-1})x_{\alpha_{n-1}}$ is contained in the image of $L(h_1)$, say it is equal to $L(h_1)y_1$ for $y_1\in L_{c_1}$. The remaining term $L(g_n)x_{\alpha_l}$ is contained in the image of $c_2\mor{h_2}e$ for $d_2\in D_1$, say it is equal to $L(h_2)y_2$ (explicitly, take $c_2 := c_{\alpha_n}$, $h_2 := g_n$ and $y_2:=x_n$). Thus there exists some $d\in D_n$ and morphisms $c_i\mor{f_i}d$ through which the two morphisms $h_i$ factor, i.e. there exists $d\mor{h}e$ such that $h\circ f_i = h_i$. But then $L_d$ contains the subspaces $L(f_1)y_1$ and $L(f_2)y_2$, so it contains their intersection, which maps to $y$ under $L(h)$. Thus $y$ is in the $\CCat$-orbits of $L_d$, and $D_n$ satisfies our assumption.
\end{proof}

In all of our applications, the indexing category $\CCat$ is closely related to the category $\FI$. In order to unify the treatment of these examples, we define the notion of a general category of $\FI$ type. To do this, we need some additional terminology.
\begin{definition}[\textbf{Weak push-out}]
A \emph{weak push-out} diagram is a pullback diagram
\squareDiagram{p}{c_1}{c_2}{d}{\tilde{f}_1}{\tilde{f}_2}{f_1}{f_2}
with the following universal property: for every other pullback diagram
\squareDiagram{p}{c_1}{c_2}{z}{\tilde{f}_1}{\tilde{f}_2}{h_1}{h_2}
there exists a unique morphism $d\mor{h} z$ that makes all the relevant diagrams commute.
We call $d$ the \emph{weak push-out object} and denote it by $c_1\coprod_p c_2$. The unique map $h$ induced from a pair of maps $c_i\mor{h_i}z$ is denoted by $h_1\coprod_p h_2$.
\end{definition}
This is similar to a usual push-out, but with ``all" commutative squares replaced by only pullback squares.
\begin{remark}
The purpose of this definition is to salvage some notion of a push-out diagram, which does not exist in a category with only injective morphisms. When starting from a category that has push-outs, such as $\Set$ and $\Vect_k$, and passing to the subcategory that includes only injective maps, we lose the push-out structure. However, weak push-outs persist, and retain most of the same function.
\end{remark}

\begin{definition}[\textbf{Categories of $\FI$ type}]
A category $\CCat$ is said to be \emph{of $\FI$ type} is it satisfies the following axioms:
\begin{itemize}
\item Every morphism is a monomorphism and every endomorphism is an isomorphism.
\item For every object $d$ there are at most a finite number of isomorphism classes of objects $c$ that satisfy $c\leq d$.
\item For every pair of objects $(c,d)$, the natural action of the group of automorphisms $\Aut_\CCat(d)$ on the set of morphisms $\Hom_\CCat(c,d)$ is transitive.
\item $\CCat$ has pullbacks and weak push-outs, i.e. for every pair of morphisms $p\mor{f_i}c_i$ there exists a weak push-out $c_1\coprod_p c_2$; and for every pair $c_i \mor{g_i} d$ there exists a pullback $c_1 \times_{d} c_2$.
\end{itemize}
We denote the group of automorphisms of an object $c$ by $G_c$.
\end{definition}

The finite number of isomorphism classes $\leq d$ implies that categories of $\FI$ type satisfy a \emph{descending chain condition}. This property will prove very useful in later treatment, and it is precisely stated as follows.
\begin{prop}[\textbf{Descending chain condition}]
Every nonempty collection $X$ of objects of $\CCat$ contains a least element, i.e. there exists some $c_0\in X$ such that if $c\in X$ satisfies $c\leq c_0$ then $c_0 \cong c$.
\end{prop}
\begin{proof}
Pick any object $d\in X$ and consider the set $X_d$ of isomorphism classes $\leq d$. Since this is a finite set, it contains a least isomorphism class, say represented by $c_0$. Now if $c\in X$ satisfies $c\leq c_0$ then by transitivity $c\leq d$, whereby we find that $c\in X_d$. But $c_0$ represents a minimal isomorphism class in $X_d$, so it follows that $c\cong c_0$.
\end{proof}

We can now formulate a criterion for the combinatorial stability of the associated intersection semilattice.
\begin{thm}[\textbf{Finite-generation implies poset stability}]
Suppose $\CCat$ is a category of $\FI$ type and that $\mathcal{A}_\bullet$ is a finitely-generated, normal $\CCat$-arrangement. Then the induced intersection semilattice $L^\mathcal{A}_\bullet$ is combinatorially stable.
\end{thm}
\begin{proof}
We have already seen in Lemma \ref{lem:normal_is_stable} that a normal $\CCat$-arrangement gives rise to a downward-stable intersection semilattice. It thus remains to show that when $\CCat$ is of $\FI$ type, the $\CCat$-poset $L^\mathcal{A}_\bullet$ is finitely-generated.

By Lemma \ref{lem:filtering_is_fg} it will suffice to show that a category of $\FI$ type is weakly filtering. For every triple of objects $p, a_1, a_2$ and a pair of maps $p \mor{\tilde{f}_i} a_i$, where $i=1,2$, we can form a weak push-out
\squareDiagram{p}{a_1}{a_2}{a_1\coprod_p a_2}{\tilde{f}_1}{\tilde{f}_2}{f_1}{f_2}
Because $\CCat$ is of $\FI$ type, the automorphisms of $a_i$ act transitively on incoming maps, therefore replacing $p\mor{\tilde{f}_i} a_i$ by some other morphism $\tilde{g}_i$ amounts to post-composing with some automorphism $\varphi_i \in G_{a_i}$. But such a replacement of $\tilde{f}_i$ with $\varphi_i\circ \tilde{f}_i$ for $i=1,2$ results in an isomorphic weak push-out object, with isomorphism given by the universal property as $\varphi_1\coprod_{p} \varphi_2$. Thus the weak push-out $a_1\coprod_p a_2$ is uniquely determined by $p$ up to isomorphism.

Let $c_1$ and $c_2$ be two objects of $\CCat$. By the definition of $\FI$ type, there exist only finitely many objects $p$ that admit maps into $c_1$, up to isomorphism. Thus by the previous argument there are only finitely many isomorphism types of weak push-outs involving $c_1$ and $c_2$. Pick representatives for these isomorphism classes $d_j = c_1\coprod_{p_j} c_2$ for $j=1,\ldots, k$. Then for every pair of morphisms $c_i\mor{g_i}e$ we can form their pullback
\squareDiagram{p}{c_1}{c_2}{e}{\tilde{f}_1}{\tilde{f}_2}{g_1}{g_2}
Now by the universal property of a weak push-out, there exists a (unique) morphism $c_1\coprod_p c_2 \mor{g}e$ that satisfies $g\circ f_i = g_i$. Find $1\leq j \leq k$ such that $d_j \cong c_1\coprod_p c_2$, then the pair $g_1$ and $g_2$ factors through $d_j$ via this isomorphism. We have thus shown that $\CCat$ is indeed weakly filtering, as every pair of morphisms $c_i\mor{f_i}e$ factors through one of the objects $d_1,\ldots,d_k$. This completes the proof.
\end{proof}
\begin{note}[\textbf{Explicit degrees of generators}] \label{rem:generating_degree_iterated_pushout}
For the purpose of finding generators explicitly in the case of categories of $\FI$ type, we trace through the construction of the sets $D_n$ from Lemma \ref{lem:filtering_is_fg}. $D_n$ is constructed inductively from $D_{n-1}$ and $D_1$ through the process of listing all possible weak push-outs. Thus if the \mbox{$\CCat$-arrangement} $\mathcal{A}_\bullet$ is generated by subspaces of $V^{c_1},\ldots, V^{c_n}$, the set $D_n$ on which all codimension $n$ generators are obtained, is the collection of objects of the form
\begin{equation}
\left(\ldots \left( c_{i_1}\coprod_{p_2} c_{i_2} \right) \coprod_{p_3} c_{i_3} \ldots \right) \coprod_{p_n} c_{i_n}
\end{equation}
\end{note}

As discussed in Step 2 on $\CCat$-posets, a combinatorially stable intersection semilattice gives rise to cohomology groups that are finitely-generated as a $\CCat$-module. We have thus proved the following. 
\begin{thm}[\textbf{Cohomology preserves finite-generation}]
If $\CCat$ is a category of $\FI$ type and $\mathcal{A}_\bullet$ is finitely-generated, normal $\CCat$-arrangement, then for every $i\geq 0$ the cohomology groups $\Ho^i(\mathcal{M}_\mathcal{A}; \Q_\ell)_\bullet$ form a finitely-generated $\CCat$-module.
\end{thm}

\subsection{Step 4 - Freeness} \label{step:freeness}
Throughout this step we assume that $\CCat$ is a category of $\FI$ type. In the context of such categories the preorder relation between objects is essentially an order:
\begin{lem} \label{lem:preorder_is_order}
Suppose $\CCat$ is a category with every morphism a monomorphism and every endomorphism an isomorphism. If $c$ and $d$ are two objects of $\CCat$ and they satisfy $c\leq d$ and $d\leq c$, then every morphism between them is an isomorphism.
\end{lem}
\begin{proof}
Let $c\mor{r}d$ and $d\mor{s}c$ be any two morphisms. We will show that $s$ has a two-sided inverse.

Consider the composition $g_c := s\circ r \in \End_\CCat(c) = G_c$. Since every endomorphism is an isomorphism, there exists an inverse $g_c^{-1}\in G_c$. We find an equality of morphisms
$$
s\circ \Id_d = s = \Id_c \circ s = \left[(s\circ r)\circ g_c^{-1}\right] \circ s = s\circ (r\circ g_c^{-1} \circ s ) 
$$
which implies that $\Id_d = (r\circ g_c^{-1}) \circ s$ as $s$ is a monomorphism. We therefore see that $r\circ g_c^{-1}$ is a two-sided inverse to $s$.
\end{proof}

We define the notion of a free $\CCat$-module through left adjoint functors. Namely, for every object $d$ there is a natural restriction functor
\begin{equation}
\Mod{\CCat} \mor{(\cdot)\rvert_{d}} \Mod{k[G_d]}
\end{equation}
that sends a $\CCat$-module $V_\bullet$ to its evaluation at $d$, i.e. $V_\bullet \mapsto V_d$. Following \cite{tD}, this functor admits a left-adjoint, which we now define.
\begin{definition} [\textbf{Induction $\CCat$-modules}] \label{def:Ind_functor}
Let $\Ind_d: \Mod{k[G_d]} \mor{} \Mod{\CCat}$ be the functor that sends a $G_d$-representation $V$ to the $\CCat$-module
\begin{equation}
\Ind_d(V)_\bullet = k[\Hom(d,\bullet)]\otimes_{G_d} V
\end{equation}
where morphisms in $\CCat$ act on these spaces naturally through their action on $\Hom(d,\bullet)$.

We call a $\CCat$-module of this form an \emph{induction module} of degree $d$.
\end{definition}
\begin{fact}
The functor $\Ind_d$ is a left-adjoint to the restriction functor at $d$. See \cite{tD} for a proof.
\end{fact}

In this context the notion of freeness is clear.
\begin{definition} [\textbf{Free $\CCat$-modules}]
A $\CCat$-module is \emph{free} if it is a direct sum of induction modules.

We define the \emph{degree} of a free module to be a least object\footnote{Least with respect to the preordering between objects, given in Definition \ref{def:ordering_objects}.} $d$ (if such exists) that is $\geq$ all the degrees of the induction modules that appear in its direct sum decomposition. If there is no such $d$ we say that the degree is $\infty$.
\end{definition}
A reader familiar with $\FI$-modules should note the free $\FI$-modules are precisely the $\FI\#$-modules introduced in \cite{CEF}.

\begin{thm}[\textbf{Free cohomology}]\label{thm:free_cohomology}
If $\mathcal{A}_\bullet$ is a continuous, normal $\CCat$-arrangement then the cohomology groups $\Ho^i(\mathcal{M}_{\mathcal{A}})$ form a free $\CCat$-module.
\end{thm}

The proof of this theorem proceeds in steps. First, we give a more useful characterization of induction $\CCat$-modules in terms of binomial sets, which we now define.
\begin{definition}[\textbf{Binomial set}] \label{def:binom_set}
Let $d$ and $e$ be two objects of $\CCat$. The group of automorphisms $G_d$ acts on $\Hom_\CCat(d,e)$ on the right by precomposition. Denote the quotient $\Hom_\CCat(d,e)/G_d$ by $\binom{e}{d}$. We will call this the \emph{binomial set, $e$ choose $d$}. If $d\mor{f}e$ is a morphism, we denote its class in $\binom{e}{d}$ by $[f]$.
\end{definition}
Note that in the case of $\CCat=\FI$, the binomial set $\binom{n}{k}$ is naturally in bijection with the collection of size $k$ subsets of $n$, hence the terminology.

\begin{lem}[\textbf{Structure of induction modules}] \label{lem:induction_form}
Suppose $M_\bullet$ is a $\CCat$-module of the form
$$
M_d = \oplus_{[f]\in\binom{d}{c}} V_{[f]}
$$
for every object $d$, and that a morphism $d_1\mor{g}d_2$ sends the factor $V_{[f]}$ to $V_{[g\circ f]}$ isomorphically. Then $M = \Ind_d(V)$ where $V := M_d$.
\end{lem}
\begin{proof}
We check that $M_\bullet$ satisfies the defining universal property of $\Ind_d(V)$. For every morphism $d_1\mor{g}d_2$ there is an isomorphism $V_{[f]} \mor{\rho_g} V_{[g\circ f]}$ for all $[f]\in\binom{d_1}{c}$, and these satisfy the compatibility condition $\rho_{h\circ g} = \rho_h \circ \rho_g$ by functoriality.

Given any $\CCat$-module $N_\bullet$ and a $G_c$-homomorphism $M_d = V \mor{\psi} N_d$, we show that there exists a unique morphism of $\CCat$-modules $M_\bullet \mor N_\bullet$ that restricts to $\psi$. For existence define a morphism ${M_d \mor{\psi_d} N_d}$ by
\begin{equation} \label{eq:induction_extension}
V_{[f]} \ni v \mapsto N(f)\circ \psi\circ \rho_f^{-1}(v)
\end{equation}
so that the diagram
\squareDiagram{V}{N_c}{V_{[f]}}{N_d}{\psi}{\rho_f}{N(f)}{\psi_e}
commutes. Here $c\mor{f}d$ is any morphism representing the class $[f]\in \binom{d}{c}$. This is well-defined since a different choice of representative $f' = f\circ g$ with $g\in G_c$ gives rise to the same map
$$
v \mapsto N(f\circ g) \circ \psi \circ \rho_{f\circ g}^{-1}(v) = N(f) \circ (N(g)\circ \psi \circ \rho_g^{-1}) \circ \rho_f^{-1}(v) = N(f)\circ \psi\circ \rho_f^{-1}(v)
$$
since $\psi$ is $G_d$-linear.

Moreover, this map is a homomorphism of $\CCat$-modules, since if $d_1\mor{h}d_2$ is any morphism then
$$
V_{[h\circ f]} \ni v \mapsto N(h\circ f)\circ\psi\circ\rho_{h\circ f}^{-1}(v) = N(h) \circ (N(f)\circ \psi \circ \rho_f^{-1})\circ \rho_h^{-1}(v) = N(h) \circ \psi_{d_1} \circ \rho_h^{-1} (v)
$$
which demonstrates that $\psi_{d_2} \circ \rho_h = N(h) \circ \psi_{d_1}$, as desired.

Lastly, this $\psi_\bullet$ is unique since every homomorphism $M_\bullet \mor{} N_\bullet$ must satisfy Relation \eqref{eq:induction_extension} by naturality.
\end{proof}

The key players in the decomposition of $\Ho^i(\mathcal{M}_{\mathcal{A}})_\bullet$ as a sum of induction \mbox{$\CCat$-modules} are primitive subspaces (see Definition \ref{def:normality}). First we show that every normal \mbox{$\CCat$-arrangement} is generated by its primitive subspaces.
\begin{lem}[\textbf{Primitive generators}]\label{lem:normal_is_primitive}
If $\mathcal{A}$ is a normal $\CCat$-arrangement then every subspace $x\in L^{\mathcal{A}}_c$ is the image of some primitive subspace.
\end{lem}
\begin{proof}
Let $x\in L_c$ be any subspace and define $X$ to be the collection of all objects $e$ for which $L_e$ contains a preimage of $x$. Explicitly, $e\in X$ when there exists some $z\in L_e$ and a morphism $e\mor{f}x$ such that $L(f)z = x$.

By the descending chain condition for categories of $\FI$ type the set $X$ contains a least object $e_0$. Choose $z_0\in L_{e_0}$ and $e_0\mor{f_0}c$ satisfying $L(f_0)z_0 = x$. We claim that $z_0$ is primitive. Indeed, if not then by definition there exists some $e_1 < e_0$ and a morphism $e_1\mor{f_1}e_0$ for which $\ker (f_1)\subseteq z_0$. Since that arrangement $\mathcal{A}$ is assumed to be normal, this implies that $z_1 = V(f_1)z_0 \in L_{e_1}$ is a preimage of $z_0$. But then $z_1$ is also a preimage of $x$, whereby we find that $e_1\in X$. This is a contradiction to the minimality of $e_0$ and thus the subspace $z_0$ must indeed be primitive.
\end{proof}
The crucial observation to make regarding primitive subspaces is that they cannot appear nontrivially in the image of any induced map of posets. For this reason they shed light onto the structure of $\CCat$, e.g. they detect isomorphisms. More generally, the following useful claim shows that two subspaces will never have equal images by way of a coincidence.
\begin{lem}[\textbf{Subspaces with equal image}]\label{lem:equal_images}
Suppose $\mathcal{A}_\bullet$ is a continuous \mbox{$\CCat$-arrangement}, $z\in L_{c_0}$ is any primitive subspace and $x\in L_{c_1}$ is any subspace. If there exists some object $d$ and morphisms $c_i\mor{f_i} d$ such that $L(f_0)z=y=L(f_1)x$, then there exists a morphism $c_0\mor{\varphi}c_1$ that sends $z$ to $x$ and satisfies $f_0= f_1\circ \varphi$. In particular, if $x$ is itself primitive, the morphism $\varphi$ is an isomorphism.

The conclusion and its proof apply even when the subspace $x\subset V^{c_1}$ is not a priori assumed to belong to the arrangement. However it follows that it must.
\end{lem}
\begin{proof}
Form the pullback of the two morphisms $f_0$ and $f_1$
\squareDiagram{p}{c_0}{c_1}{d}{r_0}{r_1}{f_0}{f_1}
and consider the corresponding diagram of vector spaces
\flippedSquareDiagram{V^p}{V^{c_0}}{V^{c_1}}{V^d}{V(r_0)}{V(r_1)}{V(f_0)}{V(f_1)}
By the continuity assumption on $\mathcal{A}_\bullet$ this is a push-out diagram of vector spaces.

If we can show that $z$ contains $\ker V(r_0) $, then since $z$ is primitive this would imply that $c_0\leq p$ and in particular $p\mor{r_0}c_0$ is an isomorphism by Lemma \ref{lem:preorder_is_order}. Using the inverse to $r_0$ we find a morphism $\varphi=r_1\circ {r_0}^{-1}: c_0\mor{}c_1$ that satisfies 
$$
f_1 \circ \varphi = f_1\circ r_1 \circ {r_0}^{-1} = f_0.
$$
In particular we have
$$
L(f_1)x = y = L(f_0)z = L(f_1) L(\varphi) z.
$$
The function $L(f_1)$ is injective, since it is defined to be $V(f_1)^{-1}$ for the surjective function $V(f_1)$. Thus we will see that $x = L(\varphi)z$, which will conclude the proof.

It remains to show that $z$ contains $\ker V(r_0)$. Since all subspaces contain the origin, we have
$$
\ker V(f_1) \subset V(f_1)^{-1}(x) = y = V(f_0)^{-1}(z)
$$ 
thus it follows that $z$ contains $V(f_0)\left(\ker V(f_1)\right)$. The claim would then follow if we can prove that there is an inclusion 
$$
\ker V(r_0)\subseteq V(f_0)\left(\ker V(f_1)\right).
$$
This follows from the universal property of $V^p$ being a push-out: $V^{c_1}$ admits a well-defined map into the quotient $V^{c_0}/V(f_0)\left( \ker V(f_1) \right)$ by first lifting to $V^d$ and then mapping into $V^{c_0}$ via $V(f_0)$, thus by the universal property there exists a map
\triangleDiagram{V^{c_0}}{V^{c_0}/V(f_0)\left( \ker V(f_1) \right)}{V^p}{}{V(r_0)}{\exists !}
that makes the diagram commute. But this implies that $\ker V(r_0)$ maps to zero through the quotient map, so $\ker V(r_0)\subseteq V(f_0)\left(\ker V(f_1)\right)$, as claimed.

Lastly, if $x$ is primitive then the same argument applied in reverse shows that $c_1\leq c_0$ as well. Thus by Lemma \ref{lem:preorder_is_order}, every morphism between $c_0$ and $c_1$ is an isomorphism.
\end{proof}

\smallskip As stated in the introduction, the general philosophy behind our approach is that representation stability phenomena are the linearized reflection of combinatorial stability. We will therefore demonstrate that freeness is already exhibited at the level of $\CCat$-sets.
\begin{lem}[\textbf{Freeness: set version}] \label{lem:free_cohomology_set_version}
Suppose $\mathcal{A}$ is a continuous, normal $\CCat$-arrangement. Then the intersection semilattice $L^\mathcal{A}$ decomposes as a disjoint union of $\CCat$-subposets corresponding to the $\CCat$-orbits of the primitive subspaces of $\mathcal{A}$. Moreover, if $\CCat(z)_\bullet$ is the orbit of the primitive subspace $z\in L_c$, then $\CCat(z)_d$ decomposes as
$$
\CCat(z)_d = \coprod_{[f]\in\binom{d}{c}} \CCat(z)_{[f]}
$$
where $\CCat(z)_{[f]}$ is the set of images of $z$ under maps induced by morphisms $c\mor{f}d$ in the equivalence class $[f]$. Lastly, every morphism $d_1\mor{g}d_2$ induces a bijection $\CCat(z)_{[f]} \mor{L(g)} \CCat(z)_{[g\circ f]}$.
\end{lem}
\begin{proof}
For every primitive subspace $z\in L_c$ we consider its $\CCat$-orbit, i.e. the $\CCat$-subposet of $L_{\mathcal{A}}$ described by
$$
\CCat(z)_d = \{ L(f)z \mid f\in \Hom_\CCat(c,d) \}.
$$
This is clearly closed under the action of $\CCat$ on $L^{\mathcal{A}}_{\bullet}$. We decompose this further as
$$
\CCat(z)_d = \bigcup_{[f]\in \binom{d}{c}} \CCat(z)_{[f]}
$$
where ${\CCat(z)_{[f]} = \{ L(f)z \mid f\in [f] \}}$. Again, it is clear that the poset map induced by $d\mor{g}e$ takes the set $\CCat(z)_{[f]}$ into $\CCat(z)_{[g\circ f]}$. We start by showing that this union is in fact disjoint. Suppose that there exists some $x\in \CCat(z)_{[f_0]}\cap \CCat(z)_{[f_1]}$. Then by Lemma \ref{lem:equal_images} it follows that there exists a morphism $c\mor{\varphi}c$ that satisfies $f_1 = f_0 \circ \varphi$. But since $\varphi\in \Hom_\CCat(c,c)=G_c$ we find that $f_1\sim f_0$.

Next we show that for every morphism $c\mor{f}d$ then induced map $\CCat(z)_{[\Id_c]}\mor{L(f)} \CCat(z)_{[f]}$ is a bijection. This will prove that all morphisms $d\mor{g}e$ indeed induce bijections $\CCat(z)_{[f]}\mor{L(g)}\CCat(z)_{[g\circ f]}$. We define the inverse map as follows. Consider the induced linear surjection $V_d \mor{V(f)}V_c$. The inverse function to $L(f)=V(f)^{-1}$ is the direct image under $V(f)$. To see that this indeed provides an inverse, note that since $V(f)$ is surjective it follows that for all $x\in \CCat(z)_{[\Id_c]}$
$$
x = V(f) \left(V(f)^{-1}x\right) = V(f) \left(L(f)x\right)
$$
and for the reverse composition, if $y\in \CCat(z)_{[f]}$ then there exists some morphism $f'\in [f]$ for which $L(f')z = y$. Since $f\sim f'$, there exists some $\varphi\in G_c$ such that $f' = f\circ \varphi$. Denote $x = L(\varphi)z$ and observe that $L(f)x=L(f')z=y$. It now follows that
$$
y = L(f)x = L(f) V(f)\left(L(f)x\right) = L(f)\left(V(f)y\right).
$$

Suppose that there exists some ${x\in \CCat(z_0)_{[f_0]}\cap \CCat(z_1)_{[f_1]}}$, i.e. there exist two morphisms $f_i$ for $i=1,2$ such that $L(f_0)z_0=x=L(f_1)z_1$. By Lemma \ref{lem:equal_images} there exists an isomorphism $c_0\mor{\varphi}c_1$ taking $z_0$ to $z_1$. This shows that the two orbits $\CCat(z_0)$ and $\CCat(z_1)$ coincide.

\smallskip Lastly, we need to show that $L_\bullet$ is a disjoint union of such $\CCat$-sets $\CCat(z)_\bullet$. Lemma \ref{lem:normal_is_primitive} asserts that every subspace $x\in L_d$ is the image of some primitive subspace, and thus it belongs to one of the $\CCat$-sets $\CCat(z)_\bullet$ described here. For their disjointness, assume that $z_i\in L_{c_i}$ for $i=1,2$ are two primitive subspaces such that there exists some $x\in \CCat(z_0)_{[f_0]}\cap \CCat(z_1)_{[f_1]}$, i.e. there exist two morphisms $f_i$ for $i=1,2$ such that $L(f_0)z_0=x=L(f_1)z_1$. By Lemma \ref{lem:equal_images} there exists an isomorphism $c_0\mor{\varphi}c_1$ taking $z_0$ to $z_1$, demonstrating that the two orbits $\CCat(z_0)$ and $\CCat(z_1)$ coincide.
\end{proof}

The set version of freeness suggests that we should consider distinct $\CCat$-orbits of primitive subspaces. This leads naturally to the following notion of equivalence.
\begin{definition}[\textbf{Equivalence of primitive subspaces}]
If $z_i\in L_{c_i}$ are primitive subspace with $i=1,2$, we say that $z_1\sim z_2$ if $\CCat(z_1)=\CCat(z_2)$. Equivalently, $z_1\sim z_2$ if there is a morphism $c_1\mor{f}c_2$ for which $L(f)z_1 = z_2$.

Denote the equivalence class of $z_1$ by $[z_1]$, and the set of all equivalence classes by $\mathcal{Z}$.
\end{definition}

We can now proceed with the final step of the proof: showing that the cohomology groups of a normal $\CCat$-arrangement are free.
\begin{proof}[Proof of Theorem \ref{thm:free_cohomology}]
Recall that Formula \eqref{eq:Bjorner_Ekedahl} states that
$$
\Ho^i\left( \mathcal{M}_\mathcal{A} \right)_c = \bigoplus_{x\in L_c} \widetilde{\Ho}_{2\cd(x)-i-2}(\Delta(L_c^{<x}))
$$
and that a morphism $c\mor{f}d$ acts on these expressions through the induced isomorphism of posets $L_c^{<x}\mor{L(f)}L_d^{<L(f)x}$. For the sake of brevity we denote the $i$-th cohomology group of $\mathcal{M}_{\mathcal{A}}^d$ by $\Ho^i_d$ and its summand $\widetilde{\Ho}_{2\cd(x)-i-2}(\Delta(L_c^{<x}))$ by $\Ho^i(x)$. In this notation, the map $f_*$ on cohomology, induced by a morphism $c\mor{f}d$, maps $\Ho^i(x)$ isomorphically onto $\Ho^i(L(f)x)$.

Let $L_\bullet = \coprod_{[z]\in \mathcal{Z}} \CCat(z)_\bullet$ be the disjoint union decomposition described in the set version of the statement, Lemma \ref{lem:free_cohomology_set_version}. Then the cohomology decomposes as a direct sum of $\CCat$-submodules:
\begin{equation}
\Ho^i_d = \bigoplus_{[z] \in \mathcal{Z}} \left(\bigoplus_{x\in \CCat(z)_d} \Ho^i(x)\right) =: \bigoplus_{[z]\in \mathcal{Z}} M^{[z]}_d
\end{equation}
We claim that for every $[z]\in \mathcal{Z}$ the corresponding direct summand $M^{[z]}_\bullet$ is an induction module, hence $\Ho^i_\bullet$ is free.

Indeed, the direct sum decomposes further as
\begin{equation}
M^{[z]}_d = \bigoplus_{x\in \CCat(z)_d} \Ho^i(x) = \bigoplus_{[f]\in \binom{d}{c}} \left( \bigoplus_{x\in \CCat(z)_{[f]}} \Ho^i(x) \right) =: \bigoplus_{[f]\in \binom{d}{c}} M^{[z]}_{[f]}
\end{equation}
Every morphism $d\mor{g}e$ takes the set $\CCat(z)_{[f]}$ bijectively onto $\CCat(z)_{[g\circ f]}$, and for every element $x\in \CCat(z)_{[f]}$ the map $\Ho^i(x)\mor{g} \Ho^i(L(g)x)$ is an isomorphism. Thus $g_*$ maps the summand $M^{[z]}_{[f]}$ isomorphically onto $M_{[g\circ f]}$. Thus according to the characterization of induction modules given in Lemma \ref{lem:induction_form} this is indeed an induction $\CCat$-module. If $z\in L_c$ is a representative primitive subspace of the class $[z]$ then $M^{[z]}_\bullet$ has degree $c$ and is generated by the $G_c$-representation
\begin{equation}
M^{[z]}_c = \bigoplus_{x\in \CCat(z)_c} \Ho^i(x) = \Ind_{\Stab(z)}^{G_c} \Ho^i(z)
\end{equation}
\end{proof}
For the purpose of keeping track of the degree of the free $\CCat$-module $\Ho^i(\mathcal{M}_\mathcal{A})_\bullet$, observe that the degrees of the induction modules that appear in its direct sum decomposition are the objects on which primitive subspaces are defined. Since only primitive subspaces $z$ of codimension $\frac{i}{2}\leq \cd(z)\leq i$ contribute to the cohomology groups, the degrees range only over objects that carry primitive subspaces with codimension in this range. Furthermore, from Note \ref{rem:generating_degree_iterated_pushout} we know that if the $\CCat$-arrangement $\mathcal{A}$ is generated in degrees $c_1,\ldots, c_n$, then all codimension-$(\leq i)$ primitive subspaces are in the image of iterated weak push-outs of at most $i$ many objects from this list. The following definition will make referring to the resulting degrees easier. 
\begin{definition}\label{def:plus_times_objects}
If $c$ and $d$ are two objects of $\CCat$, let $c+d$ denote a minimal (isomorphism class of) object that satisfies
$$
c+d \geq c \coprod_{p} d
$$
for every weak push-out of $c$ and $d$.

Similarly if $i\in \N$, let $i\times c$ denote a minimal (isomorphism class of) object that satisfies
$$
i\times c \geq c \coprod_{p_1} \ldots \coprod_{p_{i-1}} c
$$
for every $i$-fold weak push-out of $c$.
\end{definition}
Such a minimal objects exist in a category of $\FI$ type because of the descending chain condition. In fact, in all of our examples these objects are uniquely determined and can be identified explicitly: it will be given by weak coproduct $c+d = c\coprod_{\emptyset} d$ and $i\times c = c \coprod_{\emptyset}\ldots\coprod_{\emptyset} c$.

Using this notation we succinctly bound the degree of the resulting $\CCat$-modules.
\begin{cor} [\textbf{Bound on degree}]
If $\mathcal{A}_\bullet$ is a continuous, normal $\CCat$-arrangement generated in degrees $\leq c$. Then the degree of the free $\CCat$-module $\Ho^i(\mathcal{M}_{\mathcal{A}})$ is $\leq i\times c$.
\end{cor}
\begin{proof}
By the comment made in Note \ref{rem:generating_degree_iterated_pushout}, if $\mathcal{A}_\bullet$ is generated in degrees 
$$
c_1,\ldots,c_n \leq c
$$
then the $\CCat$-module $\Ho^i(\mathcal{M}_{\mathcal{A}})$ is generated in degrees given by their $i$-fold iterated weak push-outs. It is easy to verify that the definition of weak push-outs implies that for all $p\leq c_i, c_j$ we have a relation
$$
c_i \coprod_p c_j \; \leq \; c\coprod_p c \; \underset{\text{def.}}{\leq}\; c+c = 2\times c
$$
Then by induction we see that $i\times c$ is greater than all $i$-fold weak push-outs of the objects $c_1,\ldots,c_n$. In particular, all the primitive generators of $\Ho^i_\bullet$ must appear in degrees $\leq i\times c$.
\end{proof}

This concludes the proof that the cohomology groups form a free $\CCat$-module, and thus our main theorem is proved.

\section{Normalization and a criterion for normality} \label{sec:normalization}
The normality assumption in our Main Theorem \ref{thm:main_theorem} is meant to exclude cases where subspaces that could appear early in the $\CCat$-arrangement are omitted for some reason and only appear later. A normal $\CCat$-arrangement is saturated in the sense that every subspace that ``should" belong to it actually does.

We saw earlier in Lemma \ref{lem:normal_is_primitive} that a normal $\CCat$-arrangement is generated by primitive subspaces. The theorem we now state provides a converse. It also serves as an easily verifiable criterion for checking normality.
\begin{thm}[\textbf{Primitive generators imply normality}] \label{thm:primitive_is_normal}
Suppose $\CCat$ has pullbacks and $\mathcal{A}_\bullet$ is a continuous $\CCat$-arrangement. If $\mathcal{A}_\bullet$ is generated by primitive subspaces then it is a normal $\CCat$-arrangement.
\end{thm}
\begin{proof}
Suppose $Z=\left\{ z_\alpha \subset V^{c_\alpha} \right\}_{\alpha\in A}$ is a set of primitive subspaces that generates $\mathcal{A}_\bullet$. Let $d\mor{g}e$ be a morphism and let $y\in L_e$ be a subspace that contains $\ker V(g) $. We need to show that $y$ is in the image of $L(g)$.

By assumption $X$ generates the $\CCat$-arrangement, thus there exist morphisms $c_{\alpha_i}\mor{f_i}e$ where $1\leq i \leq l$ such that
$$
y = L(f_1)z_{\alpha_1} \cap \ldots \cap L(f_l)z_{\alpha_l}.
$$
Note that for every $1\leq i \leq l$ it follows that $\ker V(f_i) \subseteq L(f_i)z_{\alpha_i}$.

We prove the claim by induction on $l$. For $l=1$ the subspace $y\in L_e$ is the image of a primitive subspace $z\in L_c$ under some morphism $c\mor{f}e$. Consider the direct image $x:= V(g)y$. Since $\ker V(g) \subset y$, it follows that $V(g)^{-1}x = y$. Thus by Lemma \ref{lem:equal_images} there exists a morphism $c\mor{\varphi} d$ that takes $z$ to $x$. In particular $x\in L_d$ and $L(g)x = y$, as desired.

Now for the induction step. Assume that $y = y_1\cap y_2$ where each is the intersection of less than $l$ many images of primitive subspaces. Since $\ker V(g) \subset y_1, y_2$ the induction hypothesis implies that $y_1$ and $y_2$ are in the image of $L(g)$. But since $L(g)$ respects intersection $y$ is also in the image of $L(g)$. This completes the proof.
\end{proof}

The result stated in Theorem \ref{thm:main_theorem} does not apply to $\CCat$-arrangements that are not normal. However, even in the general case the result holds in the limit as the objects $c$ become sufficiently large. This follows from the following construction which we call \emph{normalization}.
\begin{thm}[\textbf{Normalization of a $\CCat$-arrangement}]
Suppose $\CCat$ is a category of $\FI$ type and $\mathcal{A}_\bullet$ is a continuous, finitely-generated $\CCat$-arrangement, with underlying diagram of vector spaces $V_\mathcal{A}^{\bullet}$. Then there exists a unique finitely-generated, normal $\CCat$-arrangement $\overline{\mathcal{A}}_\bullet$ also defined on $V_\mathcal{A}^\bullet$ that coincides with $\mathcal{A}_\bullet$ on a full subcategory that is upward closed and cofinal in $\CCat$. This normal $\CCat$-arrangement $\overline{\mathcal{A}}_\bullet$ will be called the \emph{normalization} of $\mathcal{A}$. 
\end{thm}
\begin{proof}
For the uniqueness statement, it suffices to show that any two normal $\CCat$-arrangements, whose underlying diagram of vector spaces is $V^\bullet$, and that coincide on a cofinal subcategory, are equal. Indeed, suppose $\mathcal{B}$ and $\mathcal{B}'$ are two such arrangements with corresponding intersection semilattices $L_\bullet$ and $L'_\bullet$ resp. Let $x\in L_c$ be any subspace and $c\leq d$ is an object such that $L_d = L'_d$. Pick a morphism $c\mor{f}d$, then $y = L(f)x\in L_d$ contains $\ker(V(f))$. Thus since $\mathcal{A}'$ is normal and $y\in L'_d$, the subspace $V(f)y = x$ belongs to $L'_c$. This shows that $L_c\subseteq L'_c$ for every object $c$. The same argument gives the opposite inclusion.

For existence suppose $K = \{x_i \subset V^{c_i}\}_{i=1}^n$ is a set of generators for $\mathcal{A}_\bullet$. We will define a new arrangement by specifying a generating set $\overline{K}$ as follows. For every $i$ let $X_i$ be the collection of objects $e$ for which there exists a morphism $e\mor{f}c_i$ with $\ker(V(f)) \subseteq x_i$. By the descending chain condition $X_i$ contains a least object $e_i$ that possesses such a morphism $e_i\mor{f}c_i$. Denote the image $V(f)x_i\subset V^{e_i}$ by $z_i$ and observe that since $\ker(V(f))\subseteq x_i$ we get an equality
$$
x_i = V(f)^{-1}V(f)x_i = V(f)^{-1}z_i.
$$

We claim that $z_i$ is primitive, i.e. it does not contain $\ker V(g)$ for any $e\mor{g}e_i$ with $e< e_i$. Indeed, if $e\mor{g}e_i$ is a morphism and $\ker V(g)\subseteq z_i$ then it follows that
$$
x_i = V(f)^{-1}z_i \supseteq V(f)^{-1}\ker V(g) = \ker \left(V(g)\circ V(f) \right) = \ker V(f\circ g)
$$
so by definition $e\in X_i$. But this is a contradiction to the minimality of $e_i$ in $X_i$.

Define a new set of generators $\overline{K}=\{z_1,\ldots,z_n\}$ and let $\overline{\mathcal{A}}_\bullet$ be the $\CCat$-arrangement generated by them: the underlying diagram of vector spaces is the same as that of $\mathcal{A}_\bullet$, and the intersection semilattice at an object $d$ is made up of all the subspaces of the form
$$
V(g_1)^{-1}(z_{i_1}) \cap \ldots \cap V(g_l)^{-1}(z_{i_l})
$$
for an $l$-tuple of morphisms $e_{i_j}\mor{g_j}d$ and $l\in \N$. It is straightforward to check that this indeed produced a $\CCat$-arrangement, and by construction it is generated by the set $\overline{K}$ of primitive subspaces. Lemma \ref{thm:primitive_is_normal} then shows that this arrangement is normal.

We claim that $\mathcal{A}$ is a subarrangement of $\overline{\mathcal{A}}$ and that the two coincide on all objects $d\geq c_1,\ldots, c_n$ (this is clearly an upward-closed and cofinal subcategory). For the first claim, note that for every $i$ there is a morphism $e_i\mor{f_i}c_i$ such that $L(f_i)z_i=x_i$. Thus $x_i\in L ^{\overline{\mathcal{A}}}_{c_i}$, and since these subspaces generate $\mathcal{A}$ we have containment $L^{\mathcal{A}}_d \subseteq L^{\overline{\mathcal{A}}}_d$ for every object $d$. Conversely, suppose $d$ admits maps from $c_1,\ldots, c_n$. It will suffice to show that every morphism $e_i\mor{g}d$ factors as $e_i \mor{f_i} c_i \mor{h} d$ for some $h$, for then the images of $z_i$ coincide with the images of $x_i$ in $L^{\overline{\mathcal{A}}}_d$. Indeed, pick any morphism $c_i\mor{h_0}d$ and consider the two morphisms $e_i\mor{h_0\circ f_i} d$ and $e_i\mor{g}d$. Since $G_d := \Aut_{\CCat}(d)$ acts transitively on incoming morphisms we can find an automorphism $\varphi\in G_d$ such that $\varphi\circ h_0\circ f_i = g$. Set $h = \varphi\circ h_0$, it satisfies $h\circ f_i = g$ as desired.
\end{proof}
\begin{note}
For concreteness' sake we reiterate that if $\mathcal{A}$ is generated in degrees $\leq c$, then a full subcategory on which $\mathcal{A}$ coincides with its normalization is made up of objects $d$ that satisfy $d\geq c$.
\end{note}

As immediate corollaries we find that the results of the Main Theorem \ref{thm:main_theorem} apply to general finitely-generated arrangements in degrees larger than those of the generating subspaces. The combinatorial version of this observation is the following.
\begin{thm}[\textbf{Limiting combinatorial stability}] \label{thm:limit_combinatorial_stability}
Suppose $\CCat$ is a category of $\FI$ type and $\mathcal{A}$ is a continuous, finitely-generated $\CCat$-arrangement, generated in degrees $\leq c$. Then the intersection semilattice of $\mathcal{A}$ coincides with a combinatorially stable $\CCat$-poset on the full subcategory of objects $d\geq c$. Namely, it coincides with the intersection semilattice of the normalization $\overline{\mathcal{A}}$.
\end{thm}
The representation theoretic version is the following.
\begin{thm}[\textbf{Limiting freeness of cohomology}]
If $\CCat$ and $\mathcal{A}$ are as in Theorem \ref{thm:limit_combinatorial_stability}, the cohomology groups $\Ho^i\left( \mathcal{M}_\mathcal{A} \right)$ coincide with a finitely-generated, free $\CCat$-module on the full subcategory of objects $d\geq c$.
\end{thm}

\section{Cohomological stability of arrangement quotients} \label{sec:cohomological_stability}
In a companion paper to this one (\cite{Ga}) we develop the theory of free and finitely-generated $\CCat$-modules, their characters, and representation stability in this context. This section will apply the theory to conclude cohomological stability results for quotients of $\CCat$-arrangements by their automorphism groups.

We will be assuming throughout that $\CCat$ is a category of $\FI$ type that is locally finite, i.e. has finite $\Hom$-sets. Denote the automorphism group of an object $c$ by $G_c$. The results that we shall use are the following. Proofs can be found in \cite{Ga}.
\begin{fact}[\textbf{Free $\CCat$-modules: Tensor products} {(\cite{Ga}[Theorem 3.7])}] \label{fact:tensor_of_free}
The tensor product of two free $\CCat$-modules of respective degrees $\leq c_1$ and $\leq c_2$ is again free of degree $\leq c_1+c_2$.
\end{fact}
\begin{fact}[\textbf{Free $\CCat$-module: Dualization} {(\cite{Ga}[Theorem 3.18])}] \label{fact:dual_of_free}
If $M_\bullet$ is a free (finitely-generated) $\CCat$-module of degree $\leq c$, then the duals $M_\bullet^*:=\Hom_\Q(M_\bullet, \Q)$ fit together to form another free (finitely generated) $\CCat$-module.
\end{fact}
\begin{fact}[\textbf{Free $\FI$-modules: Coinvariants} {(\cite{Ga}[Theorem 4.5])}] \label{fact:coinvariants_of_free}
The coinvariants of a free $\CCat$-module $M_\bullet$ of degree $\leq c$ form a $\CCat$-module $(M_G)_\bullet$. All morphisms $d\mor{f}e$ induce the same map $(M_G)_d\mor{}(M_G)_e$ which is always an injection, and is an isomorphism when $c\leq d\mor{f}e$.
\end{fact}
We also discuss the character polynomials of $\CCat$. These are defined as follows.
\begin{definition}
Suppose $M_\bullet$ is a $\CCat$-module with $M_c$ finite dimensional for every object $c$. The restriction $M_c$ is a $G_c$-representation. Denote its character by $\chi_c$. Collecting these class functions into a single object, we define the \emph{character} of the $M_\bullet$, denoted by $\chi_\bullet$, to be the class function simultaneously defined on all automorphism groups $G_c$.

A \emph{character polynomial} $P$ is a linear combination of characters of induction modules $\CCat$-modules (equivalently, free ones). Its \emph{degree}, $\deg(P)$, is a least object $d$ which is $\geq$ all degrees of the induction modules whose characters appear in $P$.
\end{definition}
Character polynomials have natural interpretations when $\CCat=\FI$, the category of finite sets and injective functions. In this case they are simply polynomials in the simultaneous class functions $X_1,X_2,\ldots$ where $X_k$ counts the number of $k$-cycles in any permutation of a finite set.
\begin{fact}[\textbf{Stabilization of inner products} {(\cite{Ga}[Corollary 4.7])}] \label{fact:inner_products}
If $P$ and $Q$ are character polynomials of degrees $p$ and $q$ respectively, then the $G_c$-inner product
$$
\langle P, Q \rangle_{G_c} := \frac{1}{|G_c|} \sum_{g\in G_c} P(g)Q(g^{-1})
$$
is independent of $c$ for all $c \geq p+q$.
\end{fact}
This is a direct consequence of the above facts regarding free $\CCat$-modules. With these facts in hand we turn back to the cohomology of $\CCat$-arrangements.

\bigskip Let $\mathcal{A}$ be a $\CCat$-arrangement. For every object $c$ the group $G_c$ acts on the variety $\mathcal{M}_\mathcal{A}^c$ and we can form the orbit space (or scheme) $\mathcal{M}_{\mathcal{A}}^c/G_c$. We denote the collection of quotient spaces by $\mathcal{M}_\mathcal{A}^\bullet/G_\bullet$ (even though they do not in fact form a functor).

Let us consider families of sheaves that arise from $\CCat$-modules via the following process. Fix some $\CCat$-module $N_\bullet$. For every object $d$ one can form the constant sheaf $\bar{N}_d$ on the space $\mathcal{M}_{\mathcal{A}}^d$. Now the pairs $(\mathcal{M}_\mathcal{A}^d, \bar{N}_d)$ fit together naturally into a $\CCat^{op}$-diagram of spaces plus a sheaf on each space. In particular one can apply the sheaf cohomology functor to this collection and get a $\CCat$-module.

Pushing the sheaf $\bar{N}_d$ forward to the quotient $\mathcal{M}_{\mathcal{A}}^d \mor{q^d}\mathcal{M}_{\mathcal{A}}^d/G_d$, the sheaf $q^d_*(\bar{N}_d)$ now admits a $G_d$-action.
\begin{definition}[\textbf{Twisted sheaf induced by a $\CCat$-module}] \label{def:twisted_sheaf}
The subsheaf of $G_d$-invariant sections $q^d_*(\bar{N}_d)^{G_d}$ will be called the twisted $N_d$-sheaf on $\mathcal{M}_{\mathcal{A}}^d$ and it will be denoted by $\widetilde{N}_d$. 
\end{definition}
\begin{remark}
When $G_d$ acts on the space $\mathcal{M}_\mathcal{A}^d$ freely, this construction yields (the sheaf analog of) the familiar Borel construction of a flat vector bundle $\mathcal{M}_{\mathcal{A}}^d \times_{G_d} N_d $. Otherwise $\widetilde{N}_d$ will be a constructible sheaf whose stalks might be smaller than $N_d$. One could check that for every point $x\in \mathcal{M}_\mathcal{A}^d$ there is an isomorphism
$$
(\widetilde{N}_d)_x \cong N_d^{\Stab_{G_d}(\tilde{x})}
$$
where $\tilde{x}\in (q^d)^{-1}(x)$ and the group $\Stab_{G_d}(\tilde{x})\subset G_d$ is the stabilizer of $\tilde{x}$. Different choices of points $\tilde{x}$ will produce different isomorphisms.
\end{remark}
We claim that when $N_\bullet$ is a free $\CCat$-module, cohomological stability holds with these systems of twisted coefficients.

\begin{thm}[\textbf{Twisted cohomological stability}]\label{thm:cohomological_stability_twisted}
Suppose $\CCat$ is a category of $\FI$-type whose automorphism groups are finite. Let $\mathcal{A}_\bullet$ be a continuous, normal $\CCat$-arrangement, generated in degree $\leq c$. Let $N_\bullet$ be a free $\CCat$-module over $\Q$ and suppose that it is generated in degree $\leq c'$. Then the sheaf cohomology groups
$$
\Ho^i(\mathcal{M}_{\mathcal{A}}^d/G_d ;\widetilde{N}_d)
$$
exhibit cohomological stability in the following sense: if $d\leq e$ then there is a well-defined injective map
$$
\Ho^i(\mathcal{M}_\mathcal{A}^d/G_d;\widetilde{N}_d) \hookrightarrow \Ho^i(\mathcal{M}_\mathcal{A}^{e}/G_{e};\widetilde{N}_e)
$$ 
and these maps become isomorphisms when $d\geq (i\times c)+c'$.

In particular, when considering the trivial $\CCat$-module $N_\bullet\equiv \Q$ (free of degree $0$), this yields a classical cohomological stability statement for $\Ho^i(X_d)$ in the range $i\times c \geq d$.

\smallskip When dealing with $\ell$-adic cohomology, one needs $N_\bullet$ to take values in continuous $\Q_\ell$-modules.
\end{thm}
\begin{proof}
For every object $d$ let $i_d$ denote the inclusion $\widetilde{N}_d = (q^d_* \bar{N}_d)^{G_d} \hookrightarrow q^d_*\bar{N}_d$. In the other direction define a transfer morphism $\widetilde{N}_d \overset{\tau_d}{\longleftarrow} q^d_*\bar{N}_d$ by $\tau_d = \frac{1}{|G_d|}\sum_{g\in G_d} g(\cdot)$. Clearly the composition $\tau_d\circ i_d$ is the identity map on $\widetilde{N}_d$ and the reverse composition $i_d\circ \tau_d$ is the projection onto the $G_d$-invariants of $q^d_*\bar{N}_d$ which are also the $G_d$-coinvariants.

Consider the induced maps on the sheaf cohomology
$$
\Ho^i(\mathcal{M}_{\mathcal{A}}^d/G_d; \widetilde{N}_d) \overset{\tau_d}{\underset{i_d}{\leftrightarrows}} \Ho^i(\mathcal{M}_{\mathcal{A}}^d/G_d; q^d_*\bar{N}_d) = \Ho^i(\mathcal{M}_{\mathcal{A}}^d;\bar{N}_d)
$$
They induce a natural isomorphism $\Ho^i(\mathcal{M}_{\mathcal{A}}^d/G_d; \widetilde{N}_d) \cong \Ho^i(\mathcal{M}_{\mathcal{A}}^d;\bar{N}_d)_{G_d}$ where the latter is the coinvariant quotient.

By the Universal Coefficients Theorem there is a natural isomorphism
\begin{equation}\label{eq:tensor_product}
\Ho^i(\mathcal{M}_{\mathcal{A}^d};\bar{N}_d) \cong \Ho^i(\mathcal{M}_{\mathcal{A}^d})\otimes N_d
\end{equation}
and as $d$ ranges over all objects this is the tensor product of two free $\CCat$-modules of respective degrees $\leq i\times c$ and $c'$. To compute the coinvariant quotients of these tensor products we use the results on free $\CCat$-modules quoted above. Using Fact \ref{fact:tensor_of_free}, the tensor product in \eqref{eq:tensor_product} is again free of degree $\leq (i\times c) + c'$, and using Fact \ref{fact:coinvariants_of_free} the coinvariant quotients stabilize in the desired sense.
\end{proof}

Often, one is interested in cohomological stability with coefficients in the various sequences of representations, e.g. the irreducibles $V_\lambda$ in the case of $S_n$. Such a sequence is a natural candidate for a homological stability statement e.g. if its characters are given by a single character polynomial. In this case, when one is interested only in the dimension of the sheaf cohomology groups a stronger stability statement can be phrased.
\begin{thm}[\textbf{Stabilization of twisted Betti numbers}]\label{thm:betti-numbers}
If $\mathcal{A}$ is as in Theorem \ref{thm:cohomological_stability_twisted} and $N_\bullet$ is any $\CCat$-module over $\Q$ whose character coincides with a character polynomial of degree $\leq d$ ($N_\bullet$ need not be free), then the dimensions of the sheaf cohomology groups
$$
\dim_k \Ho^i(\mathcal{M}_{\mathcal{A}^e}/G_e;\widetilde{N}_e)
$$
do not depend on $e$ for all $e \geq (i\times c)+d$.
\end{thm}
\begin{remark}
In Theorem \ref{thm:betti-numbers} there is no reference to the structure of $N_\bullet$ other than its character. For example, all morphisms $d\mor{f}e$ for $d<e$ might induce the zero map, and this will not be detected by the character. In particular, even though the cohomology groups eventually have the same dimension, there is no hope of finding natural isomorphisms $\Ho^i(\mathcal{M}_{\mathcal{A}}^d/G_d;\widetilde{N}_d)\mor{f_*}\Ho^i(\mathcal{M}_{\mathcal{A}}^{e}/G_{e};\widetilde{N}_{e})$
coming from the structure of $N_\bullet$.
\end{remark}
\begin{proof}
The argument in the proof of Theorem \ref{thm:cohomological_stability_twisted} above shows that
\begin{equation}
\Ho^i(\mathcal{M}_{\mathcal{A}}^d/G_d; \widetilde{N}_d) \cong \left(\Ho^i(\mathcal{M}_{\mathcal{A}^d})\otimes_k N_d\right)_{G_d}
\end{equation}
and the dimension of this coinvariant quotient is given by the $G_d$-inner product of characters
\begin{equation}
\dim_k \Ho^i(\mathcal{M}_{\mathcal{A}}^d/G_d; \widetilde{N}_d) = \langle \Ho^i(\mathcal{M}_{\mathcal{A}^d})^*, N_d \rangle_{G_d}.
\end{equation}
By Fact \ref{fact:dual_of_free}, since $\Ho^i(\mathcal{M}_{\mathcal{A}^e})$ is a free $\CCat$-module of degree $\leq i\times c$, so is its dual. Thus the character of $\Ho^i(\mathcal{M}_{\mathcal{A}^d})^*$ is given by a character polynomials $Q$ of degree $\leq i\times c$. By assumption there exists some character polynomial $P$ of degree $d$ that coincides with the character of $N_\bullet$. Therefore the above inner product of characters is given by the inner product $\langle Q, P\rangle_{G_e}$ which stabilizes for all $e\geq \deg(Q)+\deg(P^*) = (i \times c) + d$ by Fact \ref{fact:inner_products}.
\end{proof}

\section{Applications} \label{sec:applications}
In all of the following examples we consider complex varieties, i.e. we take $k=\C$. However, it should be noted that the same results hold in positive characteristic as well. Moreover, for varieties defined over $\Z$, the $\CCat$-modules we get are naturally isomorphic for any characteristic, so in this sense we might as well concentrate on the complex version of the statements. 

Consider the category $\FI$ of finite sets and injections. Every finite set is isomorphic to a unique set of the form
$$
n := \{0,\ldots,n-1\}
$$
and the endomorphisms of $n$ are the symmetric group on $n$ letters, $S_n$. Furthermore we consider the categorical power $\FI^m$ where $m\in \N$. Every object in this category is isomorphic to a unique $m$-tuple $\bar{n}:=(n^{(1)},\ldots,n^{(m)})$, and the automorphism group of $\bar{n}$ is the product of symmetric groups $S_{\bar{n}}:= S_{n^{(1)}}\times\ldots\times S_{n^{(m)}}$. The $+$ and $\times$ operations on objects coincide in this case with coordinatewise addition and multiplication. In \cite{Ga} we show that for every $m$ the category $\FI^m$ is of $\FI$ type and describe its character polynomials.
\begin{fact}[{\cite{Ga}[Theorem 6.11]}]
The $\Q$-algebra of character polynomials for $\FI^m$ is the polynomial ring
$$
\Q[X_1^{(1)},\ldots,X_1^{(m)},X_2^{(1)},\ldots, X_2^{(m)},\ldots]
$$
where $X^{(i)}_k$ is the simultaneous class function on all $S_{\bar{n}}$ given by
$$
(\sigma^{(1)},\ldots,\sigma^{(m)}) \mapsto \#\text{ of $k$-cycles in }\sigma^{(i)}
$$
whose degree is the $m$-tuple $k \bar{e}^{(i)}$.
\end{fact}

We will apply the theory developed in the paper to $\FI^m$-modules.
\begin{definition}[\textbf{The $\FI^{op}$-vector space $V^\bullet$}]
Fix some finite dimensional complex vector space $V$. We consider the $\Set^{op}$-vector space $V^\bullet: c\mapsto \Hom_{\Set}(c,V)$, i.e. $V^c$ is the $\C$-vector space of functions from $c$ to $V$ with pointwise operations.
\end{definition}
This is a continuous contravariant functor from sets to vector spaces. Moreover, since in $\Set$ every injection has a retraction, the functor $V^\bullet$ sends injections to surjective linear maps. Restricting to the subcategory $\FI$ we get an $\FI^{op}$-vector space with all induced maps surjective.
\begin{definition}[\textbf{The $\FI^m$-vector space $V^\bullet$}]
Embedding $\FI^m$ into $\Set$ naturally by considering $(A_1,\ldots,A_m) \mapsto A_1\coprod\ldots\coprod A_m$, we turn $V^\bullet$ into a $(\FI^m)^{op}$-vector space.
\end{definition}
The embedding $\FI^m\hookrightarrow \Set$ sends every morphism to an injective function, and furthermore takes pullbacks and (weak) push-outs in $\FI^m$ respectively to pullbacks and push-outs in $\Set$. In turn we find that $V^\bullet$ respects the pullbacks and weak push-outs of $\FI^m$ and sends every morphism to a surjective linear map of vector space. Thus $V^\bullet$ can serve as a continuous underlying diagram of vector spaces for $\FI^m$-arrangements.

\begin{note}
Unpacking the definition of $V^\bullet$ we see that its value at $(n^{(1)},\ldots, n^{(m)})$ is
$$
V^{n^{(1)}\coprod \ldots \coprod n^{(m)}} \cong V^{n^{(1)}}\times \ldots \times V^{n^{(m)}}.
$$
The automorphism group $S_{\bar{n}}=S_{n^{(1)}}\times\ldots\times S_{n^{(m)}}$ acts on this product through the action of every $S_{n^{(i)}}$ permuting the order of coordinates in $V^{n^{(i)}}$.
\end{note}

Using criterion \ref{thm:primitive_is_normal} for the normality of $\CCat$-arrangements we get the following result.
\begin{cor}[\textbf{Producing normal $\FI^m$-arrangements}] \label{cor:normal_FIm_arrangements}
Let $\mathcal{A}_\bullet$ be an $\FI^m$-arrangement whose underlying diagram of vector spaces is $V^\bullet$ for some $V$ and that is generated by primitive subspaces. Then $\mathcal{A}_\bullet$ is normal.
\end{cor}
The Main Theorem \ref{thm:main_theorem} then applies.
\begin{thm}[\textbf{Representation stability of $\FI^m$-arrangements}]\label{thm:application_rep_stability}
Let $\mathcal{A}_\bullet$ be any $\FI^m$-arrangement whose underlying diagram of vector spaces is $V^\bullet$ for some $V$. If $\mathcal{A}_\bullet$ is generated by primitive subspaces in degrees $\leq \bar{n}$ then
$$
\Ho^i(\mathcal{M}_\mathcal{A})_\bullet := \Ho^i(\mathcal{M}^\bullet_\mathcal{A})
$$
is a free, finitely-generated $\FI^m$-module of degree $\leq i\times \bar{n}$. In particular:
\begin{enumerate}
\item The character of $\Ho^i(\mathcal{M}_\mathcal{A})_\bullet$ is given by a character polynomial of degree $\leq i\times \bar{n}$.
\item If $P$ is any character polynomial of degree $\bar{d}$ then the character inner product $\langle P,\chi_{\Ho^i(\mathcal{M}_\mathcal{A})}\rangle_{S_{\bar{k}}}$ does not depend on $\bar{k}$ for all $\bar{k}\geq (i\times \bar{n}) + \bar{d}$. 
\item The irreducible decomposition of $\Ho^i(\mathcal{M}_\mathcal{A})_{\bar{d}}$ stabilizes in the sense of \cite{Ga} when $\bar{d}\geq 2i\times \bar{n}$.
\end{enumerate}
\end{thm}

Theorem \ref{thm:cohomological_stability_twisted} implies cohomological stability of the $S_{\bar{n}}$-quotient spaces.
\begin{thm}[\textbf{Cohomological stability of quotients}]\label{thm:application_homological_stability}
Let $\mathcal{A}_\bullet$ be as in Theorem \ref{thm:application_rep_stability} above. Consider the quotient spaces $X_{\bar{k}} = \mathcal{M}_{\mathcal{A}}^{\bar{k}}/S_{\bar{k}}$.

For any free $\FI^m$ module $N_\bullet$ of degree $\leq \bar{d}$ let $\widetilde{N}_{\bar{e}}$ is the corresponding twisted coefficient sheaf on $X_{\bar{e}}$ (see Definition \ref{def:twisted_sheaf}). Then for all $i\geq 0$ and every pair of objects $\bar{e}\leq \bar{e}'$ there is a canonical inclusion of cohomology groups
$$
\Ho^i(X_{\bar{e}}; \widetilde{N}_{\bar{e}}) \hookrightarrow \Ho^i(X_{\bar{e}'}; \widetilde{N}_{\bar{e}'})
$$
and these inclusions become isomorphisms when $\bar{e}\geq (i\times \bar{n})+\bar{d}$. In particular, for the trivial $\FI^m$-module $N_e \equiv \Q$ (which is free of degree $0$) we get classical cohomological stability with constant coefficients for all $e\geq i\times \bar{n}$, i.e. there exist canonical injections
$$
\Ho^i(X_{\bar{e}}) \hookrightarrow \Ho^i(X_{\bar{e}'})
$$
which become isomorphisms when $\bar{e}\geq i\times \bar{n}$.

Lastly, for every $i\geq 0$ and any $m$-tuple of partitions $\bar{\lambda}$ consider the $\FI^m$-module of irreducible $S_\bullet$-representations $V_{\bar{\lambda}}$ defined in \cite{Ga}. There exists a dimension $d^i_{\bar{\lambda}}$ such that for every object $\bar{e}\geq |\bar{\lambda}|+\bar{\lambda}_1, (i\times n)+|\bar{\lambda}|$ the sheaf cohomology groups satisfy
$$
\dim \Ho^i(X_{\bar{e}}; \widetilde{V}_{\bar{\lambda}(\bar{e})}) = d^i_{\bar{\lambda}}.
$$
\end{thm}

A general example to which this theory applies is the following. All later examples will be instances of this general case.
\begin{example}[\textbf{The arrangement $\mathcal{M}^\bullet_{m,k}$ and $X^\bullet_{m,k}$}]\label{example:PConf_Conf}
Let $V$ be a finite-dimensional complex vector space. Fix two natural numbers $(m,k)$ and consider the $\FI^m$-arrangement $\mathcal{A}^{m,k}_{\bullet}(V)$ whose underlying diagram of vector spaces is $V^\bullet$ and is generated by the diagonal line in $V^k\times\ldots\times V^k = V^{km}$:
\begin{equation*}
\Delta_{(k,\ldots,k)} = \left\{ \left( (z^{(1)}_1,\ldots,z^{(1)}_k),\ldots, (z^{(m)}_1,\ldots,z^{(m)}_k) \right) \mid z^{(i_1)}_{j_1}=z^{(i_2)}_{j_2} \; \forall i_1,i_2,j_1,j_2 \right\}
\end{equation*}

The preimage of $\Delta_{(k,\ldots,k)}$ under an injection $(k,\ldots,k) \mor{\bar{f}} (n^{(1)},\ldots, n^{(m)})$ is the subspace of $V^{n^{(1)}}\times\ldots\times V^{n^{(m)}}$ defined by the equations
$$
z^{(i_1)}_{f_{i_1}(j_1)} = z^{(i_2)}_{f_{i_2}(j(2))}
$$
for all $1\leq i_1,i_2\leq m$ and $1\leq j_1,j_2 \leq k$. In other words, this is the subspace in which the coordinates specified by $\bar{f}$ are all equal. As we let $\bar{f}$ range over all injections, we see that the induced arrangement on $V^{\bar{n}}$ is made up of precisely the tuples in which there exists some $z\in V$ that \underline{appears in every $V^{n^{(i)}}$ factor at least $k$ times}.

The $S_{\bar{n}}$-quotient of this arrangement (resp. its complement) is formed by forgetting the ordering of the entries in each $V^{n_i}$ factor, i.e. it is the space of unordered sets with multiplicities $U_1,\ldots,U_m\in \N[V]$ with $|U_i|=n_i$ and the intersection of all these sets contains a point with multiplicity $\geq k$ (resp. contains no point with multiplicity $\geq k$). 
\begin{definition}[\textbf{$\mathcal{M}^\bullet_{m,k}$ and $X^\bullet_{m,k}$}] \label{def:PConf_Conf}
We denote the complement of $\mathcal{A}^{m,k}_{\bar{n}}(V)$ in $V^{\bar{n}}$ by $\mathcal{M}^{\bar{n}}_{m,k}(V)$ and its $S_{\bar{n}}$-quotient by $X^{\bar{n}}_{m,k}(V)$.
\end{definition}

As stated in the introduction, these varieties are various spaces of configuration of points in $V$. We will explore this geometric aspect below in specific examples.

\begin{note}
If $k=2$ and $m=1$ then (and only then) the $S_n$ action is free. In this case the quotient map $\mathcal{M}^\bullet_{1,2} \mor{} X^\bullet_{1,2}$ is a normal $S_n$-cover. It also follows that if $N_\bullet$ is any $\FI$-module then the twisted sheaf $\widetilde{N}_n$ on $X^n_{1,2}$ is actually a vector bundle isomorphic to $\mathcal{M}^n_{1,2} \times_{S_n} N_n$.

In other cases we get a branched cover $\mathcal{M}^{\bar{n}}_{m,k}(V) \rightarrow X^{\bar{n}}_{m,k}(V)$, and the twisted coefficient sheaf $\widetilde{N}_{\bar{n}}$ on $X^{\bar{n}}_{m,k}(V)$ has different stalks above different points. This phenomenon has a natural interpretation in our case. The following example illustrates this well.
\begin{example}
Think of $\Poly^n_k(\C) := X^n_{1,k}(\C)$ as the space of degree $n$ polynomials that have no roots of multiplicity $\geq k$ (see example \ref{ex:k_equals} later), and construct the twisted coefficient sheaf corresponding to the permutation representation $S_n\curvearrowright\Q^n$.

Then the stalk over a polynomial $p\in \Poly^n_k(\C)$ can be naturally described as the $\Q$-vector space spanned freely by the distinct roots of $p$. In particular, when $p$ has multiple roots (a Zariski closed condition) this vector space will have dimension smaller than $n$.
\end{example}
\end{note}

Observe that the generating subspace $\Delta_{(k,\ldots,k)}$ is primitive, as it does not contain the kernel of any map induced by any proper injection $\bar{n}\hookrightarrow(k,\ldots,k)$. Therefore by Corollary \ref{cor:normal_FIm_arrangements}, the arrangement $\mathcal{A}^{m,k}_{\bullet}$ is a normal $\FI^m$-arrangement generated in degree $(k,\ldots,k)$. We can therefore apply our general theory and find that the cohomology groups of $\mathcal{M}^\bullet_{m,k}(V)$ satisfy representation stability, as the following theorem shows.
\begin{thm}[\textbf{Representation stability of $\mathcal{M}^\bullet_{m,k}(V)$}] \label{thm:PConf_stability}
Denote the dimension of $V$ by $r$. For every $i\geq 0$ the cohomology groups
$$
\Ho^i(\mathcal{M}_{m,k}(V))_\bullet
$$
form a free $\FI^m$-module, finitely-generated in degree $\lfloor \frac{i}{r} \rfloor\times (k,\ldots, k)$. We abbreviate this degree by $\lfloor\frac{i}{r}\rfloor k$. In particular, the following results hold:
\begin{enumerate}
\item  The character of $\Ho^i(\mathcal{M}^{\bar{n}}_{m,k}(V))$ is given by a single character polynomial of degree $\lfloor\frac{i}{r}\rfloor k$ for all $\bar{n}$.
\item If $P$ is any character polynomial of degree $\bar{d}$ the character inner product $\langle P,\chi_{\Ho^i(\mathcal{M}_{m,k})}\rangle_{S_{\bar{n}}}$ does not depend on $\bar{n}$ for all ${\bar{n}\geq \lfloor\frac{i}{r}\rfloor k + \bar{d}}$. 
\item The irreducible decomposition of $\Ho^i(\mathcal{M}^{\bar{n}}_{m,k}(V))$ stabilizes in the sense of \cite{Ga} for all $\bar{n}\geq 2\lfloor \frac{i}{r}\rfloor k$.
\end{enumerate}
\end{thm}
Furthermore, the cohomology of the quotients $X^{\bullet}_{m,k}(V)$ stabilizes:
\begin{thm}[\textbf{Cohomological stability for $X^\bullet_{m,k}$}]\label{thm:Conf_stability}
Denote $r=\dim(V)$ as above and let $N_\bullet$ be a free $\FI^m$-module of degree $\bar{d}$. For every pair $\bar{n}\leq \bar{n}'$ and every $i\geq 0$ there is a canonical inclusion of sheaf cohomology
$$
\Ho^i(X^{\bar{n}}_{m,k}(V); \widetilde{N}_{\bar{n}}) \hookrightarrow \Ho^i(X^{\bar{n}'}_{m,k}(V); \widetilde{N}_{\bar{n}'})
$$
and these inclusions become isomorphisms when $\bar{n}\geq \lfloor\frac{i}{r}\rfloor k + \bar{d}$. In particular, for the trivial $\FI^m$-module $N_{\bar{n}} \equiv \Q$ this gives classical cohomological stability with range $\bar{n} \geq \lfloor\frac{i}{r}\rfloor k$.

Moreover, for every $i\geq 0$ and an $m$-tuple of partitions $\bar{\lambda}$ consider the \mbox{$\FI^m$-module} of $S_\bullet$-irreducibles $V_{\bar{\lambda}}$ as described in \cite{Ga}. There exists a dimension $d^i_{\bar{\lambda}}$ such that for every object $\bar{n}\geq |\bar{\lambda}|+\bar{\lambda}_1, \lfloor\frac{i}{r}\rfloor k+|\bar{\lambda}|$ the sheaf cohomology groups satisfy
$$
\dim \Ho^i(X^{\bar{n}}_{m,k}(V); \widetilde{V}_{\bar{\lambda}(\bar{n})}) = d^i_{\bar{\lambda}}.
$$
\end{thm}

\begin{proof}[Proof of both Theorem \ref{thm:PConf_stability} and Theorem \ref{thm:Conf_stability}]
For the case where the dimension $r=1$, all the statements follow from our general theory: see Theorem \ref{thm:application_rep_stability} for $\mathcal{M}^\bullet$ and Theorem  \ref{thm:application_homological_stability} for $X^\bullet$.

For the general case $r\geq 1$, we recall that by Lemma \ref{lem:Deligne_bounds} a subspace $x\in L^{\mathcal{A}}$ contributes to $\Ho^i$ only if $\cd(x) \leq i$. Since the arrangement is generated by the diagonal line $\Delta_{(k,\ldots,k)}\subset V^{km}$, the subspace $x$ is the intersection of a certain number of preimages of this generating diagonal, say
$$
x = \Delta^1 \cap \ldots \cap \Delta^l
$$
is such a presentation with $l$ least. We have a sequence of proper inclusions
$$
\Delta^1 \supset \Delta^1\cap \Delta^2 \supset \ldots \supset x
$$
where at each step we find a certain number of copies of $V$. Thus the codimension of every successive pair is a positive multiple of $r$, in particular the successive codimensions are $\geq r$. This forces the codimension of $x$ to be at least $r\cdot l$. The bound $\cd(x)\leq i$ then implies that $l\leq \frac{i}{r}$. In other words: $x$ is already generated by an $\lfloor \frac{i}{r}\rfloor$-fold intersection of images of the generating diagonal. By the same argument as in Lemma \ref{lem:filtering_is_fg}, $x$ is therefore in the image of some subspace defined in degree $\lfloor \frac{i}{r}\rfloor \times (k,\ldots,k)$. This shows that the $\FI^m$-module $\Ho^i$ is generated in the stated degree. All the other statement follow from this together with the general theory presented in Theorem \ref{thm:application_rep_stability} and Theorem \ref{thm:application_homological_stability}, as in the $r=1$ case.
\end{proof}
\end{example}

\subsection{Specializing to important examples}
We now vary the parameters $m$, $k$ and $r=\dim(V)$, and specialize to concrete cases of interest. All of the following examples exhibit the stability properties described in Theorem \ref{thm:PConf_stability} and Theorem \ref{thm:Conf_stability} with the specified stable ranges.
\begin{example} [\textbf{Configurations of points in the plane and square-free polynomials}] \label{ex:braid}
Let $V = \C$ (i.e. $r=1$), $m=1$ and $k=2$. The resulting spaces form the $\FI^{op}$-space of \emph{ordered (or pure) configuration space} of distinct points in $\C$, denoted by $\PConf^\bullet(\C)$. The quotients by the action of the symmetric groups $S_n$ are the unordered configuration spaces, denoted by $\Conf^\bullet(\C)$ which, by the fundamental theorem of algebra, are naturally isomorphic to the spaces of monic square-free polynomials of degree $\bullet$, denoted by $\Poly^{\bullet}(\C)$. The isomorphism is explicitly given by sending an $n$-tuple of distinct points to the unique monic, degree $n$ polynomial which vanishes precisely at these points.

For every natural number $n$, the space $\Conf^n(\C)$ is aspherical and its fundamental group is Artin's braid group on $n$ strands, denoted by $B_n$. By forgetting the ordering $\PConf^n(\C)\rightarrow\Conf^n(\C)$ is a normal $S_n$-cover corresponding to the short exact sequence
$$
1 \mor{} P_n \mor{} B_n \mor{} S_n \mor{} 1
$$
where $P_n$ is the pure braid group. Since both spaces are aspherical, they serve as classifying spaces for $B_n$ and $P_n$ respectively, and their cohomology coincides with the group cohomology.

This sequence of spaces has been intensely studies starting with Arnol'd (\cite{Arnol'd}) and Fuks (\cite{Fu}), and more recently it served as the catalyst for the development of the theory of representation stability in \cite{CF} and later of $\FI$-modules in \cite{CEF}.

In this context Theorem \ref{thm:PConf_stability} gives a new proof for representation stability first demonstrated in \cite{CF}.
\begin{thm}[\textbf{Rep. stability for $P_n$}]
For every $i\geq 0$, the group cohomology $\Ho^i(P_n)$ forms an $\FI$-module that is finitely-generated and free of degree $\leq 2i$. In particular, there exists a character polynomial $P\in \Q[X_1,X_2,\ldots]$ of degree $2i$ which coincides with the $S_n$-character of $\Ho^i(P_n)$ for all $n$.
\end{thm}
At the level of $S_n$-quotient spaces we get cohomological stability with various systems of twisted coefficients. Every $S_n$-representation naturally becomes a \mbox{$B_n$-representation} via the natural projection $B_n\rightarrow S_n$. Thus we can consider the group cohomology of $B_n$ with coefficients in $S_n$-representations. See \cite{Ga} for the definitions and notation that are used below.
\begin{thm}[\textbf{Cohomological stability for $B_n$}]
Let $V_\bullet$ be a free $\FI$-module of degree $d$ over $\Q$. In particular, $V_n$ is an $S_n$-representation. For every $i\geq 0$ there exists canonical injections
$$
\Ho^i(B_n;V_n) \hookrightarrow \Ho^i(B_{n+1};V_{n+1})
$$
and these are isomorphisms when $n\geq 2i+d$.

Furthermore, for every partition $\lambda$ and $i\geq 0$ there exists a dimension $d^i_\lambda$ such that if $n\geq |\lambda|+\lambda_1, 2i+|\lambda|$ then 
$$
\dim \Ho^i(B_n; V_{\lambda(n)}) = d^i_\lambda.
$$
\end{thm}
\end{example}

A similar example arises by considering a vector space $V$ of higher dimension.
\begin{example}[\textbf{Configurations of points in even-dimensional Euclidean space}] \label{ex:Conf_r}
Fix any $r\geq 1$ and consider $V = \C^r \cong \R^{2r}$. By taking $m=1$ and $k=2$ in Example \ref{example:PConf_Conf} we get in degree $n$ the space $\PConf^n(V)$ of ordered configurations of $n$ distinct points in $V$. The $S_n$-quotient is the space $\Conf^n(V)$ of unordered configurations of $n$ points in $V$.

The cohomology ring of $\PConf^n(V)$ was computed by F. Cohen in \cite{Co} where it was shown to have a similar structure to that of $\PConf^n(\C)$ with all degrees multiplied by $r$. These spaces arise as a local model for ordered and unordered configurations of points in smooth complex varieties and even dimensional orientable manifolds. In \cite{To}, Totaro uses the cohomology of $\PConf^n(V)$ to compute the cohomology of the space of ordered configurations $\PConf^n(M)$ where $M$ is a smooth projective complex variety.

Specializing Theorems \ref{thm:PConf_stability} and \ref{thm:Conf_stability} to this case we get a new proof of the following results, previously proved in \cite{CEF}.
\begin{thm}[\textbf{Rep. stability of configuration space of all affine spaces}]
Let $r\geq 1$ be any dimension and $\PConf^n(\C^r)$ be the space of ordered configurations of $n$ distinct points in in $\C^r$. For every fixed $i\geq 0$ the cohomology groups $\Ho^i(\PConf^n(\C^r))$ form an $\FI$-module that is free and finitely-generated in degree $2\lfloor \frac{i}{r}\rfloor$. In particular, there exists a character polynomial $P\in \Q[X_1,X_2,\ldots]$ of degree $2\lfloor \frac{i}{r}\rfloor$ which coincides with the $S_n$-character of $\Ho^i(\PConf^n(\C^r))$ for all $n$.
\end{thm}
At the level of unordered configurations we get the following.
\begin{thm}[\textbf{Cohomological stability for configuration spaces}]
Let $V_\bullet$ be a rational free $\FI$-module of degree $d$. In particular, $V_n$ is an $S_n$-representation. We consider the cohomology with coefficients in the local system $\widetilde{V}_n$ associated with the flat vector bundle
$$
\PConf^n(\C^r) \times_{S_n} V_n \mor{} \Conf^n(\C^r).
$$
For every $i\geq 0$ there exists canonical injections
$$
\Ho^i(\Conf^n(\C^r);\widetilde{V}_n) \hookrightarrow \Ho^i(\Conf^{n+1}(\C^r);\widetilde{V}_{n+1})
$$
which are in fact isomorphisms when $n\geq 2\lfloor\frac{i}{r}\rfloor+d$.

Furthermore, for every partition $\lambda$ and $i\geq 0$ there exists a dimension $d^i_\lambda$ such that if $n\geq |\lambda|+\lambda_1, 2\lfloor\frac{i}{r}\rfloor+|\lambda|$ then 
$$
\dim \Ho^i(\Conf^n(\C^r);\widetilde{V}_{\lambda(n)}) = d^i_\lambda.
$$
\end{thm}
\end{example}

\begin{example} [\textbf{The $k$-equals arrangement}] \label{ex:k_equals}
Let $V = \C$, $m=1$ and $k\geq 2$ be arbitrary in Example \ref{example:PConf_Conf}. The resulting $\FI$-arrangement was previously called the \emph{$k$-equals arrangements} (see e.g. \cite{BW}). The complement of this arrangement parametrizes all ordered configurations of points in the plane with possible coincidences, but where no $k$ points are allowed to coincide. Taking the quotient by the action of the symmetric group we get the unordered version of such configurations. By assigning a configuration of $n$ points to the unique monic, degree $n$ polynomial that vanishes on these points (with the specified multiplicity), we get an isomorphism from the $n$-th quotient $X^n_{1,k}(\C)$ to the space $\Poly^n_k(\C)$ of monic degree $n$ polynomials that have no root of multiplicity $\geq k$. These spaces of polynomials are the complements of the natural stratification of the space of all monic, degree $n$ polynomials ($\cong\A^n$), where we filter based on the maximal multiplicity of their roots.

The intersection semilattice of the $k$-equals arrangement is usually denoted by $\Pi_{n,k}$, and is isomorphic to the lattice of partitions of $\{1,\ldots,n\}$ such that every non-singleton block has size at least $k$. These posets were studied in \cite{FNRS} relating to complexes of disconnected $k$-graphs, and by Vassiliev in connection with homotopy classification of links. 

The real version of these subspace arrangements comes up in problems of computational complexity. Consider the following problem:
\begin{quotation}
\em Given real numbers $x_1,\ldots,x_n$, decide whether at least $k$ of them are equal.
\end{quotation}
Put in other words, we are asking whether the vector $(x_1,\ldots,x_n)\in \R^n$ belongs to the $k$-equals arrangement. Bj\"{o}rner-Lov\'{a}sz show in \cite{BL} that if one tries to solve this problem using the computational model of a linear decision tree\footnote{See \cite{BL} for a definition.}, then the size and depth of this tree can be bounded from below by expressions involving the Betti numbers of (the complement of) the real $k$-equals arrangement.

The Betti numbers of the real and complex arrangements and their complements are computed in \cite{BW}. Among other things, our theory provides a new proof that the Betti numbers are polynomial in $n$. The following results were not previously stated.
\begin{thm}[\textbf{Rep. stability for the $k$-equals arrangement}]
Fix some natural number $k$, and denote the complement of the $k$-equals arrangement in $\C^n$ by $M^n_k$. These form an $\FI^{op}$-space as we let $n$ vary.

For every fixed $i\geq 0$ the cohomology groups $\Ho^i(M^n_k)$ form an $\FI$-module that is free and finitely-generated in degree $k\cdot i$. In particular, there exists a character polynomial $P\in \Q[X_1,X_2,\ldots]$ of degree $k\cdot i$ which coincides with the $S_n$-character of $\Ho^i(M^n_k)$ for all $n$.

The evaluation $P(id_n) = \dim \Ho^i$ shows that the $i$-th Betti number is given by a degree-$(k\cdot i)$ polynomial in $n$.
\end{thm}
For the quotient spaces we get the following cohomological stability result.
\begin{thm}[\textbf{Stability of polynomials with bounded root-multiplicity}]
For every natural number $k$, let $\Poly^n_k(\C)$ be the space of monic, degree $n$ polynomials whose roots have multiplicity smaller than $k$. Let $V_\bullet$ be a free $\FI$-module of degree $d$ over $\Q$, and let $\tilde{V}_n$ be the associated twisted sheaf introduced in \ref{def:twisted_sheaf}.

For every $i\geq 0$ there exists canonical injections
$$
\Ho^i(\Poly^n_k(\C);\widetilde{V}_n) \hookrightarrow \Ho^i(\Poly^{n+1}_k(\C);\widetilde{V}_{n+1})
$$
which are in fact isomorphisms when $n\geq k\cdot i + d$. In particular, for $V_n \equiv \Q$ we get classical cohomological stability for $\Poly^n_k(\C)$ in the range $n\geq k\cdot i$.

Furthermore, for every partition $\lambda$ and $i\geq 0$ there exists a dimension $d^i_\lambda$ such that if $n\geq |\lambda|+\lambda_1, k\cdot i+|\lambda|$ then 
$$
\dim \Ho^i(\Poly^n_k(\C);\widetilde{V}_{\lambda(n)}) = d^i_\lambda.
$$
\end{thm}
\end{example}

\begin{example} [\textbf{Based rational maps $\Proj^1\mor{}\Proj^{m-1}$}] \label{ex:rational_m}
Let $V=\C$, ${k=1}$ and ${m\geq 2}$ be arbitrary in Example \ref{example:PConf_Conf}. The resulting space at degree $${\bar{n}=(n^{(0)},\ldots, n^{(m-1)})}$$ consists of $m$-tuples of ordered configurations of points in the plane (with possible coincidences) whose sizes are $\bar{n}$ and who do not all have a point in common. The $S_{\bar{n}}$-quotient is the unordered version which is naturally isomorphic to the space of $m$-tuples of monic polynomials $(p_0(t),\ldots, p_{m-1}(t))$, of degrees given by $\bar{n}$, such that the gcd of the polynomials in the tuple is $1$. An equivalent description of an orbit is given by considering the algebraic function it defines
$$
[p_0(t):\ldots:p_{m-1}(t)]: \Proj^1 \mor{} \Proj^{m-1}.
$$

When restricting to the case $m=2$ and to objects of the form $(n,n)$, the quotient space is naturally isomorphic to the space of rational maps $\Proj^1\mor{}\Proj^1$ of degree $n$ that are \emph{based} in the sense that they send $\infty$ to $1$. We denote the resulting space by $\Rat^n_*(\C)$. This space is the key to understanding the space $\Rat^n(\C)$ of all degree $n$ rational maps, since there is a fibration sequence 
$$
\Rat^n_* \mor{} \Rat^n \mor{\ev_\infty} \Proj^1
$$
where the latter map is evaluation at $\infty$ and the fiber over $1 = [1:1]\in \Proj^1$ is precisely $\Rat^n_*$.

The sequence of spaces $\Rat^n_*(\C)$ was studied by Segal (see \cite{Se}), where its integral cohomological stability was demonstrated. Our techniques shows that the sequence of $S_n\times S_n$-covers (obtained by choosing orderings on the zeros and poles) satisfies representation stability rationally, and in particular rational cohomological stability follows. In fact, here we extend the rational stability result to the $2$-dimensional sequence in which we allow $n_0$ and $n_1$ to vary independently.

For the case $m>2$, one considers a similar restriction to degrees of the form $(n,\ldots,n)$, in which case we get the space of degree $n$ rational maps $\Proj^1\mor{}\Proj^{m-1}$ that are based, i.e. send $\infty\in \Proj^1$ to $[1:\ldots:1]$. We denote this space by $\Rat^n_{m*}(\C)$. As in the $m=2$ case, this space is the key to understanding the space of all degree $n$ rational maps from $\Proj^1$ to $\Proj^{m-1}$ through the fibration sequence
$$
\Rat^n_{m*}(\C) \mor{} \Rat^n_m(\C) \mor{\ev_{\infty}} \Proj^{m-1}
$$
where $\ev_{\infty}$ is the evaluation at $\infty$ function.

Specializing Theorems \ref{thm:PConf_stability} and \ref{thm:Conf_stability} to this case we get new cohomological stability results for spaces of (based) rational maps. We abbreviate $(n,\ldots,n)$ by $n$.
\begin{thm}[Cohomological stability for based rational maps]
Let $V_\bullet$ be free $\FI^m$-module of degree $\bar{d}$, and let $\widetilde{V}_n$ be the associated twisted coefficient sheaf of $\Rat^n_{m*}(\C)$ (see Definition \ref{def:twisted_sheaf}).

For every $i\geq 0$ there exists a canonical injection
$$
\Ho^i(\Rat^n_{m*}(\C);\widetilde{V}_n)\hookrightarrow \Ho^i(\Rat^{n+1}_{m*}(\C);\widetilde{V}_{n+1})
$$
and these become isomorphisms when $n\geq i+\bar{d}$.

Moreover, for every $i\geq 0$ and an $m$-tuple of partitions $\bar{\lambda}$, there exists a dimension $d^i_{\lambda}$ such that for all $n\geq |\bar{\lambda}|+\bar{\lambda}_1, i+|\bar{\lambda}|$
$$
\dim \Ho^i(\Rat^n_{m*}(\C); \widetilde{V}_{\bar{\lambda}(n)}) = d^i_{\lambda}.
$$
\end{thm}
A non-trivial example of a twisted coefficient sheaf (which is not a local system) is the sheaf whose stalks above a based rational map $f$ is the $\Q$-vector space spanned freely by the \underline{distinct} $m$-tuples $(a_0,\ldots,a_{m-1})$, where $f(a_i)$ is contained in the hyperplane $z_i=0$. This is the sheaf associated to the free, degree-$1$ $\FI^m$-module $\Ind_1(\Q)$.
\end{example}

We can also consider all of the above examples with $\C$ replaces by $\C^r$. They all satisfy representation stability (Theorem \ref{thm:PConf_stability} for the ordered version) and cohomological stability (Theorem \ref{thm:Conf_stability} for the unordered version) with improved stability ranges as $r$ grows.


\end{document}